\documentclass[12pt]{amsart}
\usepackage{latexsym, amsfonts,mathrsfs, amsmath, amssymb, amscd, epsfig}
\usepackage{mathrsfs}

\newtheorem{thm}{Theorem}
\newtheorem{defn}{Definition}
\newtheorem{prop}[thm]{Proposition}
\newtheorem{lem}[thm]{Lemma}
\newtheorem{cor}[thm]{Corollary}
\newtheorem{rmk}{Remark}

\newenvironment{pf}{{\noindent \it \bf Proof:}}{{\hfill$\Box$}\\}

\def\beq#1\eeq{\begin{equation}#1\end{equation}}
\def\balign #1 #2 \ealign{\begin{aligned} #1 #2  \end{aligned} }

\def\rhosupone{\rho_{sup}^{(1)}}
\def\rhosubone{\rho_{sub}^{(1)}}
\def\rhosuptwo{\rho_{sup}^{(2)}}
\def\rhosubtwo{\rho_{sub}^{(2)}}

\def\Esupone{E_{sup}^{(1)}}
\def\Esubone{E_{sub}^{(1)}}
\def\Esuptwo{E_{sup}^{(2)}}
\def\Esubtwo{E_{sub}^{(2)}}

\def\brhop{\bar{\rho}_+}

\def\bEp{\bar{E}_+}
\def\bup{\bar{u}_+}

\def\ep{\varepsilon}

\def\trho{\tilde\rho}
\def\tb{\tilde b}

\def\s{\sigma}

\def\tmD{\tilde{\mathscr{D}}}
\def\mL{\mathscr{L}}

\def\mA{\mathscr{A}}
\def\mC{\mathscr{C}}

\def\mF{\mathscr{F}}
\def\mG{\mathscr{G}}
\def\mD{\mathscr{D}}

\def\ms{\mathfrak{s}}
\def\bxi{{\mathbf \xi}}

\def\bX{\mathbf{X}}

\def\tx{\tilde{x}}

\def\bY{\bar{Y}}
\def\tPhi{\tilde{\Phi}}
\def\tmD{\tilde{\mathscr{D}}}
\def\hPhi{\hat{\Phi}}
\def\hmD{\hat{\mathscr{D}}}

\def\lbn{\big\|\|}
\def\rbn{\|\big\|}

\def\mBV{\mathbb{V}}

\def\mbN{\mathbb{N}}

\begin{document}

\title[ ]{Stability of Transonic Shock Solutions
for One-Dimensional Euler-Poisson Equations}

\author[ ]{Tao Luo}
\address{Department of Mathematics and Statistics, Georgetown University, USA}
\email{tl48@georgetown.edu}
\author[ ]{Jeffrey Rauch}
\author[ ]{Chunjing  Xie}
\address{Department of Mathematics, University of Michigan, 530 Church
Street, Ann Arbor, MI 48109, USA} \email{rauch@umich.edu} \email{cjxie@umich.edu}
\author[ ]{Zhouping Xin}
\address{The Institute of Mathematical Sciences and department of
mathematics, The Chinese University of Hong Kong, Hong Kong}
\email{zpxin@ims.cuhk.edu.hk}

%\date{ }\today

\begin{abstract}

In this paper, both structural and dynamical stabilities of steady
transonic shock solutions for one-dimensional Euler-Poission system
are investigated. First, a steady transonic shock solution with 
supersonic backgroumd charge is shown to be structurally stable with
respect to small  perturbations of the background charge, provided
that the electric field is positive at the shock location. Second, any steady transonic shock solution with the supersonic background charge is proved to be dynamically and exponentially stable with respect to small perturbation of the initial data, provided the electric field is not too negative at the shock location.
The proof of the first stability result relies on a monotonicity argument for the shock position and the downstream
density, and a stability analysis  for subsonic and supersonic
solutions.  The dynamical stability of the steady transonic shock for the
Euler-Poisson equations can be transformed to the global well-posedness of a free boundary
problem for a quasilinear second order equation with nonlinear
boundary conditions. The analysis for the associated linearized problem plays an essential role.

\end{abstract}

\maketitle

\section{ Introduction and Main Results}

The following system of one dimensional Euler-Poisson equations:
\begin{equation}\label{UnsteadyEP}
\begin{cases}
& \rho_t+(\rho u)_x=0, \\
& (\rho u)_t+(p(\rho)+\rho u^2)_x=\rho E, \\
& E_x=\rho-b(x),%
\end{cases}%
\end{equation}
models several physical flows including the propagation of electrons in
submicron semiconductor devices and plasmas (cf. \cite{MarkRSbook})( hydrodynamic
model), and the biological transport of ions for channel proteins (cf. \cite{Shu}). In the hydrodynamical model of semiconductor devices or plasma, $u,
\rho$ and $p$ represent the macroscopic particle velocity, electron density and
pressure, respectively, $E$ is the electric field, which is generated by the
Coulomb force of particles. $b(x)>0$ stands for the density of fixed,
positively charged background ions. The biological model describes the
transport of ions between the extracellular side and the cytoplasmic side of
the membranes(\cite{Shu}). In this case, $\rho$, $\rho u$ and $E$ are the ion
concentration, the ions translational mass, and the electric field,
respectively.

In this paper, we prove two distinct stability results for steady transonic shocks These are solutions of the following time-independent equations
\begin{equation}\label{SteadyEP}
\begin{cases}
& (\rho u)_{x}=0, \\
& (p(\rho )+\rho u^{2})_{x}=\rho E, \\
& E_{x}=\rho -b(x).%
\end{cases}
\end{equation}%
The first result concerns the stability of such solutions under perturbation of a {\it constant} background charge density. That is a purely stationary result. The second concerns the global in time stability for solutions whose initial data are close to a stationary solution.

We assume that $p$ satisfies:
\begin{equation}\label{Pressureassumption}
p(0)=p^{\prime }(0)=0,\,\,p^{\prime }(\rho )>0,\,\,p^{\prime \prime }(\rho )\geq 0,\,\text{
for}\, \rho >0,\,\,p(+\infty )=+\infty.
\end{equation}

First,
consider the boundary value problem for (\ref{SteadyEP}) in an interval $0\leq
x\leq L$ with the boundary conditions:
\begin{equation}\label{BC}
(\rho ,u,E)(0)=(\rho _{l},u_{l},E_l),\quad (\rho ,u)(L)=(\rho
_{r},u_{r}).
\end{equation}%
We assume $u_{l}>0$ and $u_{r}>0.$ By the first equation in (\ref{SteadyEP}), we know
that $\rho u(x)=constant(0\leq x\leq L)$, so the boundary data must satisfy
\begin{equation}
\rho _{l}u_{l}=\rho _{r}u_{r}.
\end{equation}
If one denotes
\begin{equation}
\rho _{l}u_{l}=\rho _{r}u_{r}=J,
\end{equation}
then $\rho u(x)=J(0\leq x\leq L)$ and the velocity is given by
\begin{equation}\label{velocity}
u=J/\rho .
\end{equation}
Thus the boundary value problem  (\ref{SteadyEP}) reduces to
\begin{equation}\label{SimpleSteadyEP}
\begin{cases}
& (p(\rho )+\frac{J^{2}}{\rho })_{x}=\rho E, \\
& E_{x}=\rho -b(x),%
\end{cases}
\end{equation}%
with the boundary conditions:
\begin{equation}\label{SimpleBC}
(\rho ,E)(0)=(\rho _{l}, E_l ),\quad \rho (L)=\rho _{r}.
\end{equation}%
Use the terminology from gas dynamics to call $c=\sqrt{p^{\prime }(\rho )}
$ the sound speed. There is a unique solution $\rho =\rho _{s}$ for the
equation
\begin{equation} \label{sonic}
p^{\prime }(\rho )=J^{2}/\rho^2,
\end{equation}%
which is the sonic state (recall that $J=\rho u$).  Later on, the flow is called supersonic if
\begin{equation}
p^{\prime }(\rho )<J^{2}/\rho^2,\ i.e.\ \rho <\rho _{s}.
\end{equation}%
Similarly, if
\begin{equation}
p^{\prime }(\rho )>J^{2}/\rho^2,\ i.e. \ \rho >\rho _{s},
\end{equation}%
then the flow is said to be subsonic.

A piecewise smooth solution $(\rho, E)$ with $\rho>0$  of (8) (or equivalently (2) with $u=\frac{J}{\rho}$) is said to be a transonic shock solution, if it is separated by a shock discontinuity, and of the form
\begin{equation*}
(\rho, E)=\left\{
\begin{aligned}
(\rho_{sup}, E_{sup})(x),\,\,\,\, 0<x<x_0,\\
(\rho_{sub}, E_{sub})(x),\,\,\,\, x_0<x<L,
\end{aligned}
\right.
\end{equation*}
satisfying the Rankine-Hugoniot conditions
\begin{equation*}
p(\rho_{sup}(x_0-))+\frac{J^2}{\rho_{sup}(x_0-)}= p(\rho_{sub}(x_0+))+\frac{J^2}{\rho_{sub}(x_0+)}, E_{sup}(x_0-)=E_{sub}(x_0+),
\end{equation*}
and is supersonic behind the shock and subsonic ahead of the shock, i.e.,
 \begin{equation}
\rho_{sup}(x_0-)<\rho_s<\rho_{sub}(x_0+),
\end{equation}
and
\begin{equation}
\rho_{sub}(x)>\rho_s, \,\, \text{for}\,\,x\in[x_0, L].
\end{equation}
(We assume that the velocity is always positive.)

Note that $(\ref{SimpleSteadyEP})_1$ is singular at the sonic state ($
p^{\prime}(\rho_s)-\frac{J^2}{\rho_s^2}=0)$ and the coefficient of $\rho_x$
changes sign for the supersonic flow and subsonic flow. This makes the
problem of determining which kind of boundary conditions should be posed to
make the boundary value problem well-posed a subtle one. In the previous
works, some pure subsonic or supersonic solutions are obtained for both
1-dimensional and multidimensional cases (cf. \cite{DeMark1d, DeMark3d, Peng}). For
a viscous approximation of transonic solutions in 2-d case for the equations
of semiconductors, see \cite{GambaMorawetz}. However, there have been only a few
results for the transonic flow. In the following, we list several results
which are closely related to the present paper. First, a boundary value
problem for (\ref{SimpleSteadyEP}) was discussed in \cite{MarkPhase} for a linear pressure
function of the form $p(\rho)=k\rho$, furthermore, the  boundary conditions read  $
\rho(0)=\rho(L)=\bar \rho$ with $\bar \rho$ being a subsonic state  and the density of the background charge  satisfied $0<b<\rho_s$. The solution in \cite{MarkPhase} may contain transonic
shock. On the other hand, since the boundary conditions and the pressure
function are  special in \cite{MarkPhase}, it is desired to consider the more
general boundary conditions with more general equation of states.
 In \cite{RosiniPhase}, a phase plane analysis is given
for system (\ref{SimpleSteadyEP}). However, no transonic shock solutions  were constructed in
\cite{RosiniPhase}. A transonic solution which may contain transonic shocks was
constructed by Gamba (cf. \cite{Gamba1d}) by using a vanishing viscosity
limit method. However, the solutions as the limit of vanishing viscosity may
contain boundary layers. Therefore, the question of well-posedness of the
boundary value problem for the inviscid problem can not be answered by the
vanishing viscosity method. Moreover, the structure of the solutions constructed by the
vanishing viscosity method in \cite{Gamba1d} is shown to be of  bounded
total variation and possibly contain more than one transonic shock.

A thorough study  of the transonic shock solutions for one-dimensional  Euler-Poisson equations with a constant doping profile  $b(x)=b_0$ was given in \cite{LuoXin} for both cases when $0<b_0<\rho_s$ and $b_0>\rho_s$. In the present paper, we concentrate on the case when $0<b_0<\rho_s$. In \cite{LuoXin}, when $0<b_0<\rho_s$, and if the boundary data $(\rho_l, E_l)$, $\rho_r$ and the interval length $L$ satisfy some conditions (see \cite{LuoXin} for details), then the problem (\ref{SimpleSteadyEP}) and (\ref{SimpleBC}) admits a unique transonic shock solution.  The stability of the transonic solutions obtained in \cite{LuoXin} when $b$ is a  small perturbation of given background charge $b_0=b_0(x)$,  is analyzed in the following theorem.
\begin{thm}\label{Thmexistence} Let $J>0$ be a constant, and let $b_0$  satisfy 
\begin{equation}\label{condb}
0<\min_{x\in  [0, L]} b_0(x)\leq \max_{x\in  [0, L]} b_0(x) < \rho_s
\end{equation}
and $(\rho_l, E_l)$ be a supersonic state ($0<\rho_l<\rho_s)$, $\rho_r$ be a subsonic state ($\rho_r>\rho_s$). If the boundary value problem (\ref{SimpleSteadyEP}) and (\ref{SimpleBC}) admits a unique transonic shock solution $(\rho^{(0)}, E^{(0)})$ for the case when $b(x)=b_0(x)$ ($x\in [0, L]$) with a single transonic shock located at
 $x=x_0\in (0, L)$ satisfying
 \begin{equation}\label{positiveE}
 E^{(0)}(x_0+)=E^{(0)}(x_0-)>0,
 \end{equation}
   then  there exists $\epsilon_0>0$ such that if
\begin{equation}\label{sizeofb}
\|b-b_0\|_{C^{0}[0, L]}=\epsilon \leq \epsilon_0,
\end{equation}
 then the boundary problem (\ref{SimpleSteadyEP}) and (\ref{SimpleBC}) admits  a unique transonic shock solution $(\tilde{\rho}, \tilde{E})$ with a single transonic shock locating at some $\tilde{x}_0 \in [x_0-C\epsilon, x_0+C\epsilon]$ for some constant $C>0$.
\end{thm}

\begin{rmk}
When $b_0(x) \equiv const$, it should be noted that there are a large class of boundary data which ensure the existence and uniqueness of the transonic shock solutions satisfying the assumptions in Theorem \ref{Thmexistence} , see \cite{LuoXin}.
\end{rmk}

Second, we would like to investigate the dynamical stability of the steady transonic shock solutions.
For a given  function $b(x)$ satisfying  $0<b=b(x)<\rho_s$ for $x\in [0, L]$,  and a constant  $\bar J>0$,  let
\begin{equation}\label{stead}
\begin{aligned}
(\bar\rho,\bar{u},\bar{E})(x)=\left\{
\begin{array}{ll}
(\bar\rho_-,\bar J/\bar{\rho}_-,\bar{E}_-)(x),\quad\text{if}\,\,0<x<x_0,\\
(\bar\rho_+,\bar J/\bar{\rho}_+,\bar{E}_+)(x),\quad\text{if}\,\,x_0<x<L\\
\end{array}
\right.
\end{aligned}
\end{equation}
be a steady transonic shock solution of (\ref{SteadyEP}) satisfying the boundary conditions
\begin{equation}\label{boundary1}
(\bar \rho, \bar E)(0)=(\rho_l, E_l),\,\,\,\,  \bar \rho(L)=\rho_r.
\end{equation}
Precisely, we assume that $(\bar\rho,\bar{u},\bar{E})$ is supersonic
for $0\le x\le x_0$, and subsonic for $x_0<x\le L$, i.e.,
\begin{equation}\label{supsub1}\begin{cases}& p'(\bar \rho_-)(x)<\bar J^2/\bar{\rho}_-^2(x), \ {\rm for~} 0\le x\le x_0,\\
& p'(\bar \rho_+)(x)>\bar J^2/\bar{\rho}_+^2(x), \ {\rm for~} x_0\le x\le L. \end{cases}\end{equation}
At $x=x_0$, the Rankine-Hogoniot conditions are satisfied:
\begin{equation}\label{rh1}
\left(p(\bar \rho_-)+\frac{\bar{J}^2}{\bar{\rho}_-}\right)(x_0)=\left(p(\bar \rho_+)+\frac{\bar{J}^2}{\bar{\rho}_+}\right)(x_0), \bar E_-(x_0)=\bar E_+(x_0).
\end{equation}
Finally, we assume the solution  is away from vacuum:
\begin{equation}\label{vacuum}
\inf_{x\in [0, L]}\bar \rho(x)>0.
\end{equation}
It is easy to see that we can extend $(\bar\rho_-, \bar E_-)$ to be a smooth supersonic solution of (\ref{SteadyEP}) on $[0, x_0+\delta]$ for some $\delta>0$, which coincides with $(\bar\rho_-, \bar E_-)$ on $[0, x_0]$. Later on, we still use $(\bar\rho_-, \bar E_-)$ to denote this extended solution. Similarly, we will denote $(\brhop, \bEp)$ to be a subsonic solution of (\ref{SteadyEP}) on $[ x_0-\delta, L]$ for some $\delta>0$, which coincides with $(\brhop, \bEp)$ in (\ref{stead}) on $[x_0, L]$.

We consider the initial boundary value problem of system (\ref{UnsteadyEP})  with the initial data
\begin{equation}\label{initial}
(\rho, u, E)(0, x)=(\rho_0,u_0,E_0)(x), \end{equation} and the
boundary condition
\begin{equation}\label{bc2}
(\rho, u, E)(t, 0)=(\rho_l,\frac{\bar J}{\rho_l},E_l),\,\,\,\,
\rho(t, L)=\rho_r, \end{equation} where $\rho_l, E_l$ and $\rho_r$
are the same as in (\ref{boundary1}).

We assume that the initial data are of the form
\begin{equation}\label{in1}
(\rho_0,u_0)(x)=\begin{cases}
&(\rho_{0-},u_{0-})(x),\quad\text{if}\,\,0<x<\tilde{x}_0,\\
&(\rho_{0+},u_{0+})(x),\quad\text{if}\,\,\tilde{x}_0<x<L,
\end{cases}
\end{equation}
and
\begin{equation}\label{in2}
E_0(x)=E_l +\int_0^x (\rho_0(s)-b(s))ds
\end{equation}
which is a small perturbation of $(\bar\rho, \bar u, \bar E)$ in the sense that
\begin{equation}\label{AssID}
\begin{aligned}
&|x_0-\tilde{x}_0|+\|(\rho_{0+},u_{0+}) -(\bar\rho_+,\bar u_+)\|_{H^{k+2}([\check{x}_0, L])} \\
&+ \|(\rho_{0-},u_{0-}) -(\bar\rho_-,\bar u_-)\|_{H^{k+2}([0,\hat{x}_0])}<\varepsilon,
\end{aligned}
\end{equation}
for some small $\varepsilon>0$, and some integer $k\ge 15$, where $\check{x}_0=\min\{x_0, \tilde{x}_0\}$ and $\hat{x}_0=\max\{x_0, \tilde{x}_0\}$. Moreover, $(\rho_0, u_0, E_0)$ is assumed to satisfy the Rankine-Hogoniot conditions as $x=\tilde{x}_0$,
\begin{equation}\label{rh2}
\begin{aligned}
&\left((p(\rho_{0+})+\rho_{0+}u_{0+}^2-(p(\rho_{0-})+\rho_{0-}u_{0-}^2)\right)
\cdot(\rho_{0+}-\rho_{0-})(\tilde{x}_0)\\
=&(\rho_{0+}u_{0+}-\rho_{0-}u_{0-})^2(\tilde{x}_0).
\end{aligned}
\end{equation}

Before stating our dynamical stability result, we give the definition for piecewise smooth entropy solution for Euler-Poisson equations. 

\begin{defn}
If $(\rho_-, u_-, E_-)$ and  $(\rho_+, u_+, E_+)$ are $C^1$ smooth solutions of Euler-Poisson equations (\ref{UnsteadyEP}) in the regions $\{(t, x)|t\geq 0,  0\leq x\leq s(t)\}$ and $\{(t, x)|t\geq 0, s(t)\leq x\leq L\}$, respectively. Then 
\begin{equation}
(\rho,u,E)=\left\{
\begin{array}{ll}
(\rho_-,u_-,E_-),\quad\text{if}\quad 0<x<s(t),\\
(\rho_+,u_+,E_+),\quad\text{if}\quad s(t)<x<L
\end{array}
\right.
\end{equation}
is said to be a piecewise smooth entropy solution of (\ref{UnsteadyEP}) at $x=s(t)$ if $(\rho, u, E)$ satisfies Rankine-Hugoniot conditions
\begin{equation}\label{rh3}
\left\{
\begin{aligned}
&(p(\rho)+\rho u^2)(t,s(t)+)-(p(\rho)+\rho u^2)(t, s(t)-)\\
 &\,\,\,\,=(\rho u(t, s(t)+)-\rho u(t, s(t)-))\dot s(t)\\
&\rho u(t, s(t)+)-\rho u(t, s(t)-)=(\rho (t, s(t)+)-\rho(t, s(t)-))\dot{s}(t), \\
& E(s(t)+, t)=E(s(t)-, t),
\end{aligned}
\right.
\end{equation}
and Lax geometric entropy condition 
\begin{equation}
\begin{aligned}
&(u-\sqrt {p'(\rho)})(t, s(t)-)>\dot s(t)>(u-\sqrt {p'(\rho)})(t, s(t)+),\\
& (u+\sqrt {p'(\rho)})(t, s(t)+)>\dot s(t).
\end{aligned}
\end{equation}
\end{defn}
The dynamical stability theorem in this paper is as follows.

\begin{thm}\label{Thmdystability}
Let $(\bar\rho,\bar{u},\bar{E})$ be a steady transonic shock solution to system (\ref{UnsteadyEP}) satisfying (\ref{stead}), (\ref{boundary1}), (\ref{supsub1}), (\ref{rh1}), and (\ref{vacuum}). Moreover, there exists a $\delta>0$ ($\delta$ depends on $(\bar\rho, \bar{u}, \bar{E})$) such that
\begin{equation}\label{stabilitycondition}
\bar E_-(x_0)=\bar E_+(x_0)>-\delta. \\
\end{equation}
Then there exists an $\ep_0>0$ such that for any $\ep\leq \ep_0$,
if the initial data $(\rho_0, u_0, E_0)$ satisfy (\ref{in2}),  (\ref{AssID}), (\ref{rh2})
and the $k+2$-th order compatibility conditions at $x=0$, $x=x_0$ and $x=L$,
then the initial boundary value problem (\ref{UnsteadyEP}),   (\ref{initial}) and (\ref{bc2}) admits  a unique piecewise smooth entropy solution $(\rho, u, E)(x, t)$ for $(t, x)\in [0, \infty)\times [0, L]$ containing a single transonic shock $x=s(t)$ ($0<s(t)<L$)  with $s(0)=\tilde{x}_0$.

Furthermore, there exist $T_0>0$ and $\lambda>0$ such that
\begin{equation}
\begin{aligned}
(\rho_-,u_-,E_-)(t, x)=(\bar\rho_-,\bar u_-,\bar E_-)(x), \quad \text{for}\,\, 0\le x<s(t),
\end{aligned}
\end{equation}
for $t>T_0$ and
\begin{equation}
\begin{aligned}
\|(\rho_+,u_+,)(\cdot,t)-(\bar\rho_+,\bar u_+)(\cdot)\|_{W^{k-7,\infty}(s(t), L)} + \|E_+(\cdot,t)-\bar E_+(\cdot)\|_{W^{k-6,\infty}(s(t), L)} \leq C\varepsilon e^{-\lambda t},
\end{aligned}
\end{equation}
\begin{equation}
\begin{aligned}
\sum_{m=0}^{k-6}|\partial^{m}_t(s(t)-x_0)|\leq C\varepsilon e^{-\lambda t},
\end{aligned}
\end{equation}
for $t\ge 0$, where  $(\bar{\rho}_{\pm}, \bar{u}_{\pm},\bar{E}_{\pm})$ are the solutions of the Euler-Poisson equations in the associated regions.
\end{thm}

The condition (\ref{stabilitycondition}) is used to prove the exponential dynamical stability of the steady transonic shock solutions. When this condition is violated, we have the following linear instability results for some special cases.

\begin{thm}\label{Thminstability}
There exist $L>0$ and a linearly unstable transonic shock solution $(\bar{\rho}, \bar{u}, \bar{E})$ satisfying  (\ref{stead}), (\ref{boundary1}), (\ref{supsub1}), (\ref{rh1}), and (\ref{vacuum}) and
\begin{equation}
\bar{E}_-(x_0)=\bar{E}_+(x_0)<-C
\end{equation}
for some positive constant $C$.
\end{thm}

Begin with several remarks concerning Theorem \ref{Thmexistence} , Theorem \ref{Thmdystability}, and Theorem \ref{Thminstability}.
\begin{rmk}
In both Theorem \ref{Thmexistence} and Theorem \ref{Thmdystability}, the results are also true if we impose small perturbations for the boundary conditions (\ref{SimpleBC}).
\end{rmk}
\begin{rmk}
It follows from the results in \cite{LuoXin} and Theorem \ref{Thmexistence}, that the background transonic shock solution in Theorem \ref{Thmdystability} does exist.  Moreover, we do not assume that $b(x)$ is a small perturbation of a constant in both Theorem  \ref{Thmexistence} and Theorem \ref{Thmdystability}, where it may have large variation.
\end{rmk}
\begin{rmk}In \cite{RosiniStability}, the {\it  local-in-time stability} of transonic shock solutions for
the Cauchy problem of (\ref{UnsteadyEP}) is considered by assuming the existence of
steady transonic shocks. Here, we prove the global-in -time exponential stability for the initial boundary value problem.  \end{rmk}

\begin{rmk}
The compatibility conditions for the initial boundary value problems for hyperbolic equations were discussed in detail in \cite{RauchM, Majda, Metivier}.
\end{rmk}

\begin{rmk}
In Theorem \ref{Thmdystability}, the regularity assumption is not optimal. By adapting the method in \cite{Metivier}, less  regularity assumption than that in (\ref{AssID}) will be enough . However, our proof only involves the elementary weighted energy estimates rather than paradifferential calculus.
\end{rmk}

As far as the Euler-Poisson equations are concerned, there have been  many studies on the large time behavior of solutions, see \cite{Huang, LuoNX} and references therein. However, all these studies are for the Euler-Poisson equations with relaxations, where the relaxation  has very strong dissipative effects.

It is interesting to compare these results with the transonic solutions of a
quasi-one-dimensional gas flow through a nozzle. The time-dependent equations for the quasi-one-dimensional
isentropic nozzle flow are
\begin{equation}\label{Q1DNozzle}
\begin{cases}
& \rho_t+(\rho u)_x=-\frac{A^{\prime}(x)}{A(x)}\rho u, \\
& (\rho u)_t+ (\rho u^2 +p(\rho))_x=-\frac{A^{\prime}(x)}{A(x)}\rho u^2,%
\end{cases}%
\end{equation}
where $\rho, u$ and $p$ denote respectively the density, velocity
and pressure, $A(x)$ is the cross-sectional area of the nozzle. In
\cite{EmbidGM}, steady state solutions for (\ref{Q1DNozzle})
containing transonic shocks were constructed for the boundary value
problem in the interval [0, 1] with the boundary conditions $ (\rho,
u)(0)=(\rho_l, u_l)$ and $(\rho, u)(1)=(\rho_r, u_r)$, where
$(\rho_l, u_l)$ is  supersonic and $(\rho_r, u_r) $ is subsonic and
satisfies $\rho_lu_l=\rho_ru_r$. The general wave patterns for flows
in quasi-one-dimensional nozzles were studied in \cite{LiuARMA}. The
stability of transonic shocks for system (\ref{Q1DNozzle}) was
studied in \cite{LiuTP}. It was shown that the transonic shock
solution is stable if $A'>0$ in \cite{LiuTP} via a modified Glimm
scheme by introducing some steady states in the building blocks. For
the studies on the solutions of general hyperbolic conservation laws
with moving source, see \cite{LiuARMA, LiuJMP, Lien, Ha, HaYang} and
references therein. For piecewise smooth initial data, when the
nozzle is a straight slowly increasing sectional nozzle, Xin and Yin
obtained the dynamical stability of steady weak transonic shock
solutions, see \cite{XinYinJDE}.

\begin{rmk}
The idea for the proof of Theorem \ref{Thmdystability} can also be applied for the stability of transonic shock solutions with arbitrary amplitude for quasi-one-dimensional nozzle flows with rapidly increasing nozzle walls, see \cite{RXX}.
\end{rmk}

It would be  interesting to extend the results in this paper to the
multi-dimensional case, as those for the gas dynamics, see \cite{CF1, XinYinCPAM, LXY1, LXY2} and references therein. An effort in this direction was made in \cite{GambaMorawetz} for a viscous
approximation of transonic solutions in 2-d case for the equations of
semiconductors. However, zero viscosity limit in \cite{GambaMorawetz} remains an open problem.

The proof of persistence of transonic shocks under perturbation of $b(x)$ depends on a monotonic dependence of the shock location as a function of downstream density proved in \S \ref{secmono}. 
The a priori estimates  for supersonic and subsonic solutions and existence of transonic solutions are presented in  \S \ref{secest}.
The proof  that for initial data close to steady transonic shock the solution decays exponentially to the shock depends on the fact that 
\begin{itemize}
\item perturbations to the left as swept upstream by the supersonic hypothesis, and importantly that
\item perturbations to the right decay because of absorption at the shock.
\end{itemize}
The proof of the latter property is not straightforward. In addition to the usual technical difficulties from the quasilinear structure there is fundamental difficulty that the problem involves a free boundary (shock) on the left of the subsonic region. The key is to prove decay for the linearized problem. After a nontrivial transformation that the linearized problem resembles a Klein-Gordon equation (\S \ref{secform}). For that equation we construct a nontrivial energy functional which is decreasing due to dissipative boundary conditions (\S \ref{seclinear}).  Decay to the linearized problem is proved by employing an argument of Rauch-Taylor \cite{Rauch}. Then the quasilinear technical argument comes in \S \ref{secuap}.  The linear instability of the transonic shocks is investigated in \S \ref{secinsta} when the condition (\ref{stabilitycondition}) is violated.
Finally, we have  an  appendix for the Rauch-Taylor type estimate \cite{Rauch} for the linearized problem for unsteady Euler-Poisson equations.

%\newpage
\section{Monotone Relation between the Shock Position and the Downstream
Density}\label{secmono}

In this section, we derive a monotone relation between the shock position and the downstream density (pressure)
for the steady transonic shock solutions for the Euler-Poisson equations.

For any supersonic state $(\rho, E)$ satisfying $\rho<\rho_s$, one can connect it to a unique subsonic state $(\ms(\rho),E)$ via a transonic shock, where $\ms(\rho)$ is determined by the entropy condition and  the Rankine-Hugoniot condition:
\begin{equation}\label{RHcondition}
\ms(\rho)>\rho_s\,\,\,\,\text{and}\,\,\,\,p(\ms(\rho))+\frac{J^2}{\ms(\rho)}=p(\rho)+\frac{J^2}{\rho}.
\end{equation}
This yields
\begin{equation}
\frac{d\ms}{d\rho}(\rho)=\frac{p'(\rho)-\frac{J^2}{\rho^2}}{p'(\ms(\rho))-\frac{J^2}{\ms^2(\rho)}}.
\end{equation}
Note that  smooth solutions of (\ref{SimpleSteadyEP}) satisfy
\begin{equation}
\left\{
\begin{aligned}
& \rho_x=\frac{\rho E}{p'(\rho)-\frac{J^2}{\rho^2}},\\
& E_x=\rho -b.
\end{aligned}
\right.
\end{equation}
Hence the supersonic solution $(\rho(x), E(x))$ of (\ref{SimpleSteadyEP}) satisfies
\begin{equation}\label{SubsonicDerivative}
\frac{d \ms\Big(\rho(x)\Big)}{dx}=\frac{\rho E}{p'(\ms(\rho))-\frac{J^2}{\ms^2(\rho)}}.
\end{equation}

The monotone relation between the shock position and the downstream density is given by the following lemma.
\begin{lem}\label{lemmonotone}
Let $(\rho^{(1)}, E^{(1)})$ and $(\rho^{(2)}, E^{(2)})$ be two transonic shock solutions of (\ref{SimpleSteadyEP}), and $(\rho^{(i)}, E^{(i)}) (i=1,2)$ are defined as follows
\begin{equation}
(\rho^{(i)}, E^{(i)})=\left\{
\begin{aligned}
(\rho^{(i)}_{sup}, E^{(i)}_{sup}), \,\, \text{for } \,\, 0<x<x_i,\\
(\rho^{(i)}_{sub}, E^{(i)}_{sub}), \,\, \text{for } \,\, x_i<x<L,
\end{aligned}
\right.
\end{equation}
where
\begin{equation}
\rho^{(i)}_{sup}<\rho_s<\rho^{(i)}_{sub}\,\,\ \,\text{for}\,\,\,\, i=1,2.
\end{equation}
Moreover, they satisfy the same upstream boundary conditions,
\begin{equation}
\rho^{(1)}(0)= \rho^{(2)}(0)=\rho_l,\,\,
E^{(1)}(0)=E^{(2)}(0)=E_l.
\end{equation}

If $b<\rho_s$, $x_1<x_2$ and $\Esuptwo(x_1)>0$, then
\begin{equation}
\rho^{(1)}(L)>\rho^{(2)}(L).
\end{equation}

\end{lem}

\begin{pf}
Since $(\rhosupone, \Esupone)$ and $(\rhosuptwo, \Esuptwo)$ satisfy the same ODE system and initial values,
\begin{equation}
(\rhosupone, \Esupone)=(\rhosuptwo, \Esuptwo) \,\,\,\,\text{for}\,\,\,\, x\in [0,x_1].
\end{equation}
For $x_1\le x\le x_2$, define a function $E_{\alpha}$ as follows
\begin{equation}
\left\{
\begin{aligned}
&\frac{d E_{\alpha}}{dx}=\ms(\rhosuptwo)-b,\quad {\rm for}~ x_1\le x\le x_2,\\
&E_{\alpha}(x_1)=\Esubone(x_1)=\Esupone(x_1)=\Esuptwo(x_1).
\end{aligned}
\right.
\end{equation}
Note that $\rhosuptwo<\rho_s<\ms(\rhosuptwo)$, therefore,
\begin{equation}
\Esuptwo(x)<E_{\alpha}(x)\,\,\,\,\text{for}\,\, x\in(x_1,x_2].
\end{equation}
Since $E_{\alpha}(x_1)>0$ and $b<\rho_s$, one has $E_{\alpha}(x)>0$ for $x\in (x_1,x_2]$.
For $x\in(x_1,x_2)$, the equation (\ref{SubsonicDerivative}) for $\rho_{sup}^{(2)}$ becomes
\begin{equation}
\left\{
\begin{aligned}
&\frac{d \ms(\rho_{sup}^{(2)})}{dx}=\frac{\rho_{sup}^{(2)} E_{sup}^{(2)}}{p'(\ms(\rho_{sup}^{(2)}))-\frac{J^2}{\ms^2(\rho_{sup}^{(2)})}},\\
&\ms(\rho_{sup}^{(2)})(x_1)=\rho_{sub}^{(1)}(x_1).
\end{aligned}
\right.
\end{equation}
Thus
\begin{equation}
\frac{d \ms(\rhosuptwo)}{dx}< \frac{\rhosuptwo E_{\alpha}(x)}{p'(\ms(\rhosuptwo))-\frac{J^2}{\ms^2(\rhosuptwo)}}.
\end{equation}
It follows from the fact $E_{\alpha}(x)>0$ that, on $[x_1,x_2]$, $(\ms(\rhosuptwo), E_{\alpha})$ satisfies
\begin{equation}
\left\{
\begin{aligned}
&\frac{d \ms(\rhosuptwo)}{dx}< \frac{\ms(\rhosuptwo) E_{\alpha}}{p'(\ms(\rhosuptwo))-\frac{J^2}{\ms^2(\rhosuptwo)}},\\
&\frac{d E_{\alpha}}{dx}=\ms(\rhosuptwo)-b,\\
& \ms(\rhosuptwo)(x_1)=\rhosubone(x_1),\,\,E_{\alpha}(x_1)=\Esubone(x_1).
\end{aligned}
\right.
\end{equation}
Note that $(\rhosubone, \Esubone)$ satisfies  the equations
\begin{equation}
\left\{
\begin{aligned}
&\frac{d \rhosubone}{dx}= \frac{\rhosubone \Esubone}{p'(\rhosubone)-\frac{J^2}{(\rhosubone)^2}},\\
&\frac{d \Esubone}{dx}=\rhosubone-b,\\
& \rhosubone(x_1)=\rhosubone(x_1),\,\,\Esubone(x_1)=\Esubone(x_1).
\end{aligned}
\right.
\end{equation}
The comparison principles for ODE systems (\cite{Pao}), yields that
\begin{equation}
\ms(\rhosuptwo)(x_2)< \rhosubone(x_2),\,\,\,\, E_{\alpha}(x_2)< \Esubone(x_2).
\end{equation}
Since $\Esuptwo(x_2)<E_{\alpha}(x_2)$, we have
\begin{equation}
\Esubtwo(x_2)=\Esuptwo(x_2)< E_{\alpha}(x_2)\leq \Esubone(x_2).
\end{equation}
Note that $(\rhosubone, \Esubone)$ and $(\rhosubtwo,\Esubtwo)$ solve the same ODE system on $[x_2, L]$,  by the comparison principle for ODEs again,   one has
\begin{equation}
\rhosubone(L)> \rhosubtwo(L) \,\,\,\,\text{and}\,\,\,\, \Esubone(L)> \Esubtwo(L).
\end{equation}
This finishes the proof of the Lemma.
\end{pf}

\section{A Priori Estimates and Existence of Steady Transonic Shock Solutions}\label{secest}

In this section, we prove a priori estimates for supersonic and subsonic flows via the multiplier method, which yield the existence of supersonic, subsonic, and transonic shock solutions.

A smooth solution of (\ref{SimpleSteadyEP}), $\rho$,
satisfies the following second order equation:
\begin{equation}\label{SOsupereq}
(f(\rho)\rho')'=\rho-b,
\end{equation}
where
\begin{equation}
f(\rho)=\frac{1}{\rho}\left(p'(\rho)-\frac{J^2}{\rho^2}\right).
\end{equation}

Suppose that $(\rho_0, E_0)$ is a supersonic or subsonic solution of
steady Euler-Possion equations with background charge $b_0$ and with initial data $(\rho_{I}, E_{I})$, namely,
\begin{equation}\label{bgsupereq}
\left\{
\begin{aligned}
&f(\rho)\rho'=E,\\
&E'=\rho-b_0\\
&\rho(a)=\rho_{I},\,\,\,\, E(a)=E_{I}.
\end{aligned}
\right.
\end{equation}

The following lemma gives some stability estimates for both the supersonic and the subsonic solutions of  (\ref{SOsupereq}), which are small perturbations of the solutions of (\ref{bgsupereq}).
\begin{lem}\label{Lemweightestimate}
For any interval $[a, l]\subseteq [0, L]$, let  $(\rho_0, E_0)$ be a supersonic or subsonic solution  of the problem (\ref{bgsupereq}). Then
there exists $\epsilon>0$ such that if
\begin{equation}\label{AssumptionbIC}
\|b(x)-b_0\|_{C^0[a, l]}+|\tilde{\rho}_I|+|\tilde{E}_I|<\epsilon,
\end{equation}
then there exists a unique supersonic or subsonic solution $(\rho, E)(x)$ for $x\in [a, l]$ of the problem
(\ref{SimpleSteadyEP})  with initial conditions
\begin{equation}\label{ODEIC}
\rho(a)=\rho_I+\trho_I,\,\,\,\, E(a)=E_I+\tilde{E}_I,
\end{equation}
Moreover, $(\rho, E)$ satisfies the following estimate
\begin{equation}\label{xxx}
\|\rho-\rho_0\|_{C^{1}[a, l]}< C e^{\alpha L}\epsilon,
\end{equation}
for some constants $C>0$ and $\alpha>0$.
\end{lem}
\begin{pf}
We only give the proof for the case when $\rho_0$ is subsonic on $[a, l]$ (the case when when $\rho_0$ is supersonic is completely similar). When  $\rho_0$ is subsonic on $[a, l]$, there are constants $c_0>0$ and $C_1>0$ such that
\begin{equation}\label{add1} f(\rho_0)(x)>c_0, \quad {~\rm for~} x\in [a, l], \end{equation}
and \begin{equation}\label{add2} \rho_s<\rho_0(x)<C_1, \ |\rho_0'(x)|\le C_1, \quad {~\rm for~} x\in [a, l].
 \end{equation}
 We may assume a priorily that
 \begin{equation}\label{add3} f(\rho)(x)>c_0/2, \quad {~\rm for~} x\in [a, l], \end{equation}
and \begin{equation}\label{add4} \rho_s<\rho(x)<2C_1, \ |\rho'(x)|\le 2C_1, \quad {~\rm for~} x\in [a, l].
 \end{equation}
 Once we obtain the estimate (\ref{xxx}), under this a priori assumption, the lemma can be proved by using the standard local existence theory of ODE systems and continuation argument.

Let us denote
\begin{equation}
\trho=\rho-\rho_0, \,\,\,\tb=b-b_0.
\end{equation}
Then,
\begin{equation}\label{deviationODE}
(f(\rho)\trho'+F(\rho_0,\trho)\rho_0'\trho)'-\trho=-\tb,
\end{equation}
where
\begin{equation}
F(\rho_0,\trho)=\int_0^1 f'(\rho_0+\theta\trho)d\theta.
\end{equation}

Choosing  a multiplier $K(x)=e^{-\lambda(x-a)}$ for some constant $\lambda>0$ and
multiplying both sides of (\ref{deviationODE}) by $K(\tilde{\rho}'+\tilde{\rho})$,  after some straightforward calculations, one has
\begin{equation}
\begin{aligned}
&\int_a^l-\tilde{b}K(\tilde{\rho}'+\tilde{\rho})dx=\int_a^l\{(f(\rho)\trho'+F(\rho_0,\trho)\rho_0'\trho)'-\trho\}K(\trho'+\trho)dx\\
= &\int_a^l Kf(\rho)\left(\frac{(\trho')^2}{2}\right)'  +K[(f(\rho))'+F(\rho_0, \trho)\rho_0'](\trho')^2+K\trho(f(\rho)\trho')' dx\\
&+\int_a^l\left (K(F(\rho_0,\trho)\rho_0')' +KF(\rho_0,\trho)\rho_0'
 -K\right)\trho\trho' dx\\
 & +\int_a^l \left(K(F(\rho_0,\trho)\rho_0')'-K\right) \trho^2 dx
\end{aligned}
\end{equation}
Integration by parts and the fact  that $K'(x)=-\lambda K(x)$ and $K''(x)=\lambda^2 K(x)$ yield
\begin{equation}
\begin{aligned}
&\int_a^l-\tilde{b}K(x)(\tilde{\rho}'+\tilde{\rho})dx\\
=&\int_a^l K(x)\left (\left(\frac{\lambda}{2}-1\right) f(\rho)+\frac{1}{2}(f(\rho))' +F(\rho_0, \trho)\rho_0'\right)(\trho')^2dx\\
&+\int_a^l K(x)\left(\frac{\lambda^2}{2}f -\frac{\lambda}{2}f'+(F\rho_0')'-1)\right)\trho^2 dx\\
&+\int_a^l K(x)\left ( (F(\rho_0,\trho)\rho_0')' + F(\rho_0,\trho)\rho_0'
 -1\right)\trho\trho' dx\\
&+K(l)\left(f \frac{(\trho')^2}{2}+f\trho\trho'+\lambda f\frac{\trho^2}{2}\right)\Big|_{x=l}-(f \frac{(\trho')^2}{2}+f\trho\trho'+\lambda f\frac{\trho^2}{2})\Big|_{x=a}.
\end{aligned}
\end{equation}
In view of (\ref{add3}) and (\ref{add4}), we can choose $\lambda$ sufficient large so that
\begin{equation}\label{add5}
\begin{aligned}
&\int_a^l\lambda e^{-\lambda (x-a)}(\lambda \trho^2 + (\trho')^2)dx+ e^{-\lambda (l-a)}\left(\trho'^2(l) +\lambda \trho^2(l)^2\right)\\
\leq & C\int_a^l e^{-\lambda (x-a)}\tb^2 dx+C(\lambda |\trho_I|^2+|\tilde{E_I}|^2).
\end{aligned}
\end{equation}
Here and in the following,  $C$ denotes a generic positive constant $C$.

Therefore, if (\ref{AssumptionbIC}) holds for some $\epsilon>0$, then
\begin{equation}
\int_a^l e^{-\lambda (x-a)}\tb^2 dx +C(\lambda |\trho_I|^2+|\tilde{E_I}|^2) \leq C\epsilon^2.
\end{equation}
This, together with (\ref{add5}), implies
\begin{equation}
\int_a^l \trho^2 +(\trho')^2 dx \leq C e^{\lambda (l-a)}\epsilon^2\le C e^{\lambda L}\epsilon^2 ,
\end{equation}
which in turn implies
\begin{equation}
\|\trho\|_{C^{0}[a, l]}\leq  C e^{\alpha L}\epsilon,
\end{equation}
 where $\alpha=\lambda/2$.
Using (\ref{deviationODE}), we obtain that
\begin{equation}
\|\trho'\|_{C^0[a, l]}\leq Ce^{\alpha L}\epsilon.
\end{equation}
\end{pf}

We are now in position to prove give  Theorem \ref{Thmexistence}.

By the assumption on the unperturbed transonic shock solution $(\rho^{(0)}, E^{(0)})$ for the case when $b(x)=b_0(x)$ ($x\in [0, L]$),  there exists a constant $\delta>0$ satisfying $ [x_0-\delta, x_0+\delta]\in (0, L)$, such that the initial value problem for the  ODE system:
\begin{equation}\label{0}\begin{cases}
&(p(\rho)+\frac{J^2}{\rho})_x=\rho E, \quad E_x=\rho-b_0, \\
& (\rho, E)|_{x=0}=(\rho_l, E_l)
\end{cases}\end{equation}
admits a unique smooth solution $(\rho^{\ell}, E^{\ell})(x)$ on the interval $x\in [0, x_0+\delta]$ satisfying
$0<\rho^{\ell}(x)<\rho_s$  for $x\in [0, x_0+\delta]$ and
 \begin{equation}\label{newconstruction}  E^{\ell}(x)>0,\qquad  {\rm for}~ x_0-\delta\le x\le x_0+\delta,\end{equation} where $x_0$ is the shock location for the unperturbed background transonic shock solution  $(\rho^{(0)}, E^{(0)})$ for the case when $b(x)=b_0$ ($x\in [0, L]$). Moreover, by the uniqueness for the initial value problems of ODEs, we have
 \begin{equation}\label{7}(\rho^{\ell}, E^{\ell}(x)=(\rho^{(0)}, E^{(0)})(x), {~\rm for}\ x\in [0, x_0).\end{equation}

 Now we define two functions:  $(\rho_1^{r}, E_1^{r})(x)$  for $x\in [x_0-\delta, L]$ and $(\rho_2^{r}, E_2^{r})(x)$
 for  $x\in [x_0+\delta, L]$ as follows.
 Let $(\rho_1^{r}, E_1^{r})(x)$ be the solution of the initial value problem of the following ODE,
 \begin{equation}\label{4}\begin{cases}
&(p(\rho)+\frac{J^2}{\rho})_x=\rho E, \quad E_x=\rho-b_0, {~\rm for}\ x\in (x_0-\delta, L] \\
& (\rho, E)(x_0-\delta)=(\ms(\rho^{\ell} (x_0-\delta)), E^{\ell}(x_0-\delta))
\end{cases}\end{equation}
where $\ms$ is the function defined in (\ref{RHcondition}).
It follows from (\ref{Lemweightestimate}) that the problem (\ref{4})
admits a unique smooth subsonic solution $(\rho_1^{r}, E_1^{r})(x)$ on the interval $x\in [x_0-\delta, L]$ satisfying
$\rho_1^{r}(x)>\rho_s$ and $ E_1^{r}(x)(x)>E^{\ell}(x_0-\delta)> 0$ for $x\in (x_0-\delta, L]$.

Let $(\rho_2^{r}, E_2^{r})(x)$ be the solution of the initial value problem of the following ODE,
\begin{equation}\label{4'}\begin{cases}
&(p(\rho)+\frac{J^2}{\rho})_x=\rho E, \quad E_x=\rho-b_0, {~\rm for}\quad  x\in (x_0+\delta, L], \\
& (\rho, E)(x_0+\delta)=(\ms(\rho^{\ell} (x_0+\delta)), E^{\ell}(x_0+\delta)),
\end{cases}\end{equation}
where $\ms$ is the function defined in (\ref{RHcondition}).
Again, following  from (\ref{Lemweightestimate})   that the problem (\ref{4'})
admits a unique smooth subsonic solution $(\rho_2^{r}, E_2^{r})(x)$ on the interval $x\in [x_0+\delta, L]$ satisfying
$\rho_2^{r}(x)>\rho_s$ and $ E_2^{r}(x)(x)>E^{\ell}(x_0+\delta)> 0$ for $x\in (x_0+\delta, L]$.

By the monotonicity of the shock location w.r.t. the exit density (Lemma \ref{lemmonotone}), we have
\begin{equation}\label{8}
\rho_2^{r}(L)<\rho_r<\rho_1^{r}(L).
\end{equation}

Now, we define two transonic solutions for the case when $b$  is a small perturbation of $b_0$  as
follows: Let $x_1=x_0-\delta$ and $x_2=x_0+\delta$ and
\begin{equation}\label{9}
(\hat{\rho}^{(i)},\hat{E}^{(i)})(x)=\begin{cases} &(\hat\rho_i^{\ell}, \hat E_i^{\ell})(x), \quad, {\rm for}~ 0\le x< x_i,\\
& (\hat\rho_i^{r}, \hat E_i^{r})(x) \quad, {\rm for}~  x_i<x\le L,\end{cases}\end{equation}
for $i=1, 2$, where
$(\hat \rho_i^{\ell}, \hat E_i^{\ell})(x)$ is the solution of the following problem:
\begin{equation}\label{0}\begin{cases}
&(p(\rho)+\frac{J^2}{\rho})_x=\rho E, \quad E_x=\rho-b(x), \qquad 0<x<x_i\\
& (\rho, E)|_{x=0}=(\rho_l, E_l)
\end{cases}\end{equation} and
$(\hat\rho_i^{r}, \hat E_i^{r})$ is the solution for the following problem
\begin{equation}\label{10}\begin{cases}
&(p(\rho)+\frac{J^2}{\rho})_x=\rho E, \quad E_x=\rho-b(x), {~\rm for}\ x\in (x_i, L] \\
& (\rho, E)(x_i)=(\ms(\hat\rho_i^{\ell} (x_i-)), \hat E_i^{\ell}(x_i))
\end{cases}\end{equation}
where $\ms$ the function defined in (\ref{RHcondition}).

It follows from Lemma \ref{Lemweightestimate} that, if $\|b-b_0\|_{C^0[0, L]}=\epsilon$ is sufficiently small, then
$(\hat\rho_i^{\ell}, \hat E_i^{\ell})$ and $(\hat\rho_i^{r}, \hat E_i^{r})$ are well-defined and satisfying:
\begin{equation}
\begin{aligned}
&\sup_{x\in [0, x_i)}|(\hat\rho_i^{\ell}, \hat E_i^{\ell})-(\rho_i^{\ell}, E^{\ell})|\le C\epsilon,\notag\\
&\sup_{x\in ( x_i, L]}|(\hat\rho_i^{r}, \hat E_i^{r})-(\rho_i^{r}, E^{r})|\le C\epsilon,
\end{aligned}
\end{equation}
for $i=1, 2$ , for some constant $C>0$.  In particular,
\begin{equation}
|(\hat\rho_i^{r}(L)-\rho_i^{r}(L)|\le C\epsilon,
\end{equation}
for $i=1, 2$. This, together with (\ref{8}), implies that
\begin{equation}
\hat\rho_2^{r}(L)<\rho_r<\hat\rho_1^{r}(L),
\end{equation}
 if $\epsilon$ is sufficiently small.
This shows that the transonic shock solution of the boundary problem (\ref{SimpleSteadyEP}) and (\ref{SimpleBC}) admits  a unique transonic shock solution $(\tilde{\rho}, \tilde{E})$ with a single transonic shock located at some $\tilde{x}_0 \in (x_1, x_2)$ by a monotonicity argument as follows:  For $a\in [x_1, x_2]$, we define a function $g(a)=\rho(L)$ where $\rho$ is a transonic shock solution of the system (\ref{SimpleSteadyEP}) with $\rho(0)=\rho_l$ and $E(0)=E_l$ and its shock is at $a$. By Lemma \ref{Lemweightestimate},  we know that $g(a)$ is a continuous on $[x_1, x_2]$.  Moreover, Lemma \ref{lemmonotone} implies that $g(a)$ is  strictly decreasing on $[x_1, x_2]$. Finally, the stability estimate,  Lemma {\ref{Lemweightestimate}, yields  that $\tilde{x}_0 \in [x_0-C\epsilon, x_0+C\epsilon]$ for some constant $C>0$.

{\it Proof of the uniqueness of transonic shock solutions}. Suppose, besides the transonic shock solution we constructed above,  there is another transonic shock solution for the perturbed boundary value problem (\ref{SimpleSteadyEP}) and (\ref{SimpleBC}) for the variable $b(x)$, say, $(\rho_1, E_1)$ , which contains a single transonic shock at $x=\tilde{x}_1$.  Because of the monotonicity of the shock location w.r.t the downstream density, $\tilde{x}_1$ must lie out of the region $(x_0-\delta, x_0+\delta)$.
Let $\bar{g}$ be the function defined on $[0, L]$ as  $\bar{g}(a)=\bar{\rho}(L)$ where $\bar{\rho}$ is a transonic shock solution of the system (\ref{SimpleSteadyEP}) with $\rho(0)=\rho_l$ and $E(0)=E_l$ and its shock locates at $a$ for the case when $b(x)=b_0$. By Lemma \ref{Lemweightestimate}, $\bar{g}$ is a continuous function, so
\begin{equation}
S=\{y|y=\bar{g}(a) , \,\, a \in [0,L]\backslash (x_0-\delta, x_0+\delta)\}
\end{equation}
is a closed set. Since we assume that the transonic shock solution $(\bar\rho_1, \bar E_1)$ to the problem (8) and (9) for the case when $b(x)=b_0$ is unique,  $\rho_r\notin S$. Hence, there is an $\epsilon_1>0$ such that
\begin{equation}
\inf_{y\in S}|\rho_r-y|>\epsilon_1.
\end{equation}
In particular,
\begin{equation}\label{exitdis}
|\rho_r-\bar{g}(\tilde{x}_1)|>\epsilon_1.
\end{equation}
On the other hand, by Lemma \ref{Lemweightestimate}, we have
\begin{equation}
|\rho_r-\bar{g}(\tilde{x}_1)|=
|\rho_r-\bar{\rho}_1(L)|=|\rho_1(L)-\bar{\rho}_1(L)|<\epsilon_1/2.
\end{equation}
if $b(x)$ satisfies (\ref{sizeofb}) for some $\epsilon_0>0$.
This contradicts (\ref{exitdis}). Hence, there is only one transonic shock solution for the problem (\ref{SimpleSteadyEP}) and (\ref{SimpleBC}).

%\newpage

\section{Dynamical Stability of Transonic Shock Solutions}\label{secdyna}

In this section, we investigate the dynamical stability of transonic
shock solutions for the Euler-Poisson equations
(\ref{UnsteadyEP}).

\subsection{Formulation of the Problem}\label{secform}
Let $(\bar{\rho},\bar{u},\bar{E})$ be a steady transonic shock solution of the form (\ref{stead}) satisfying (\ref{stabilitycondition}). Suppose that the initial data $(\rho_0,u_0,E_0)$   satisfies (\ref{AssID})
and the $k+2$-th order compatibility conditions.
It follows from the argument in \cite{LY}  that there exists a piecewise smooth solution containing a single shock $x=s(t)$ (with $(s(0)=\tilde{x}_0$) satisfying the Rankine-Hogoniot conditions and Lax geometric
 shock condition (\ref{rh3}), of the Euler-Poisson equations on $[0, \bar T]$ for some $\bar T>0$, which can be written as
\begin{equation}\label{ustsolution}
(\rho,u,E)(x, t)=\left\{
\begin{array}{ll}
(\rho_-,u_-,E_-),\,\,\text{if}\,\,0<x<s(t),\\
(\rho_+,u_+,E_+),\,\,\text{if}\,\,s(t)<x<L.
\end{array}
\right.
\end{equation}
Note that, when $t>T_0$ for some $T_0>0$, $(\rho_-, u_-, E_-)$ will depend only on the boundary conditions at $x=0$. Moreover, when $\varepsilon$ is small, by the standard lifespan argument, we have $T_0< \bar T$
(cf. \cite{LY}).
 Therefore,
\begin{equation}\label{ste}
(\rho_-, u_-, E_-)=(\bar{\rho}_-, \bar{u}_-, \bar{E}_-)\,\, \text{for}\,\, t>T_0.
\end{equation}

In the following, for simplicity of the presentation, we may
assume $T_0=0$ without loss of generality. We want to extend the
local-in-time solution to all $t>0$. In view of (\ref{ste}), we only
need to obtain uniform estimates in the region $x>s(t), t>0$.
For this purpose, we will formulate an initial boundary value
problem in this region. First, the Rankine-Hugoniot conditions for
(\ref{ustsolution})  read
\begin{equation}\label{USRH}
[\rho u]=[\rho]s'(t),\,\, [\rho u^2+p]=[\rho u]s'(t),
\end{equation}
where $[f]=f(s(t)+,t)-f(s(t)-,t)$, so
\begin{equation*}
[p+\rho u^2]\cdot[\rho]=[\rho u]^2.
\end{equation*}
More precisely, one has
\begin{equation*}
\begin{aligned}
&\left(p(\rho_+)(t,s(t))+\frac{J_+^2}{\rho_+}(t,s(t))-p(\rho_-)(t,s(t))
-\frac{J_-^2}{\rho_-}(t,s(t))\right)\cdot(\rho_+-\rho_-)\\
=&(J_+(t,s(t))-J_-(t,s(t)))^2,
\end{aligned}
\end{equation*}
where $J=\rho u$.  In view of (\ref{ste}),
\begin{equation}\label{JJ}
J(s(t)-, t)=\bar J.\end{equation}
Hence,
\begin{equation*}
\begin{aligned}
&\left(p(\rho_+)(t,s(t))+\frac{J_+^2(t,s(t))}{\rho_+(s(t),t)}-p(\bar\rho_+)(s(t))
-\frac{\bar J^2}{\bar\rho_+}(s(t))\right.\\
&\left.+p(\bar\rho_+)(s(t))+\frac{\bar J^2}{\bar\rho_+}(s(t))
-p(\bar\rho_-)(s(t))-\frac{\bar
J^2}{\bar\rho_-}(s(t))\right)\cdot(\rho_+(t,s(t))-\rho_-(s(t)))\\
=&(J_+(t,s(t))-\bar J)^2.
\end{aligned}
\end{equation*}
It follows from the Rankine-Hugoniot conditions for (\ref{rh1}) and Taylor expansions that
\begin{equation*}
\begin{aligned}
&\left\{p'(\bar\rho_+)(s(t))(\rho_+(t,s(t))-\bar\rho_+(s(t)))-\frac{\bar J^2}{\bar\rho_+^2}(s(t))\cdot(\rho_+(t,s(t))-\bar\rho_+(s(t)))\right.\\
&+\frac{2\bar J}{\bar\rho_+}(s(t))\cdot(J_+(t,s(t))-\bar J(s(t)))+\partial_x(p(\bar\rho_+)+\frac{\bar J^2}{\bar\rho_+})(x_0)\cdot(s(t)-x_0)\\
&\left.-\partial_x(p(\bar\rho_-)+\frac{\bar
J_-^2}{\bar\rho_-})(x_0)\cdot(s(t)-x_0)+R_1\right\}
(\bar\rho_+(x_0)-\bar\rho_-(x_0)+R_2)\\
=&(J_+(t,s(t))-\bar J(s(t)))^2,
\end{aligned}
\end{equation*}
where
\begin{equation*}
R_1=O(1)((\rho_+-\brhop)^2+(J_+-\bar J)^2+(s(t)-x_0)^2)\, .
\end{equation*}
Later on, we always use $R_1$ to denote those quadratic terms with
different $O(1)$ coefficients. It may change from line to line. Thus, by implicit function theorem, we have
\begin{equation}\label{DiffJshock}
(J_+-\bar J)(t,s(t))=\mA_1((\rho_+-\brhop)(t,s(t)), s(t)-x_0)
\end{equation}
where $\mA_1$ regarded as a function of two variables satisfies
\begin{equation*}
\begin{aligned}
\mA_1(0,0)=0,\,\, \,\, \frac{\partial \mA_1}{\partial
(\rho_+-\brhop)}=-\frac{p'(\bar\rho_+)-\frac{\bar J^2}{\bar\rho_+^2} }{2\bar J/\bar\rho_+}(x_0),\,\,\,\,
\frac{\partial \mA_1}{\partial (s-x_0)}=-
\frac{(\bar\rho_+-\bar\rho_-)\bar E_+}{2\bar J/\brhop}(x_0).
\end{aligned}
\end{equation*}
Substituting (\ref{DiffJshock}) into the first equation in (\ref{USRH}) yields
\begin{equation}\label{shockspeedeq}
s'(t)=\mA_2(\rho_+-\brhop,s(t)-x_0)
\end{equation}
where $\mA_2$ satisfies $\mA_2(0,0)=0$ and
\begin{equation*}
\frac{\partial \mA_2}{\partial(\rho_+-\brhop)}=-\frac{p'(\bar\rho_+)-\bar J^2/\rho_+^2}{2\bar
u_+(\bar\rho_+-\bar\rho_-)}(x_0),\,\, \,\,\frac{\partial \mA_2}{\partial(s(t)-x_0)}=-\frac{\bar
E_+}{2\bar u_+}(x_0).
\end{equation*}

It follows from  the third equation in
(\ref{UnsteadyEP}) that
\begin{equation*}
\begin{aligned}
E_+(x, t)=E_l+\int_0^{s(t)}(\rho_--b)(y)dy+\int_{s(t)}^x(\rho_+-b)(y)dy,
\end{aligned}
\end{equation*}
for $s(t)<x\le L$.
Applying the first equation in (\ref{UnsteadyEP})  and the Rankine
-Hugoniot conditions (\ref{USRH}), one has
\begin{equation*}
\begin{aligned}
\partial_tE_+=-\rho_+u_+(t, x)+\rho_lu_l=-J_+(t, x)+\bar J.
\end{aligned}
\end{equation*}

%\newpage

Set
\begin{equation}
\begin{aligned}
Y=E_+(x,t)-\bar E_+(x).
\end{aligned}
\end{equation}
Then
\begin{equation}
Y_t=\bar J-J_+,\,\,\,\, Y_x=\rho_+-\brhop.
\end{equation}
Therefore, it follows from the second equation in the Euler-Poisson
equations (\ref{UnsteadyEP}) that
\begin{equation*}
\begin{aligned}
(E_+-\bar E_+)_{tt}+\partial_x\left(p(\bar\rho_+)+\frac{\bar J^2}{\bar\rho_+}-p(\rho_+)-\frac{J_+^2}{\rho_+}\right)+\rho_+
E_+-\bar\rho_+\bEp=0.
\end{aligned}
\end{equation*}
Thus,
\begin{equation}\label{3120}
\begin{aligned}
\partial_{tt}Y+\partial_x\left(p(\bar\rho_+)+\frac{\bar J^2}{\bar\rho_+}-p(\bar\rho_++Y_x)-\frac{(\bar J-Y_t)^2}{\bar\rho_++Y_x}\right)+\bar E_+\partial_xY +\bar\rho_+Y+YY_x=0.
\end{aligned}
\end{equation}
One can write this equation in the following way. With $\bxi=(\xi_0, \xi_1) =(t, x)$ and
$$\partial_{i}=\frac{\partial}{\partial_i}, \partial_{ij}=\frac{\partial^2}{\partial_i\partial_j}$$
for $i, j=0, 1$,
\begin{equation}\label{EQY}
\sum_{0\le i, j\le 1}\bar{a}_{ij}(x,Y_t,Y_x)\partial_{ij}Y+\sum_{0\le i\le 1}\bar{b}_i(x,Y_t,Y_x)\partial_i Y+\bar{c}(x,Y_t, Y_x)Y=0,
\end{equation}
where $\bar a_{ij}$, $\bar b_i$ and $\bar c$ are smooth functions of their arguments, and satisfy
\begin{equation}
\begin{aligned}
\mL_0 Y &= \sum_{0\le i, j\le 1}\bar{a}_{ij}(x,0,0)\partial_{ij}Y+\sum_{0\le i\le 1}\bar{b}_i(x,0,0)\partial_i Y+\bar{c}(x,0, 0)Y\\
& = \partial_{tt}Y-\partial_x((p'(\bar\rho_+)-\frac{\bar J^2}{\bar\rho_+^2})\partial_xY)+\partial_x(\frac{2\bar J}{\bar\rho_+}\partial_tY)+\bar E_+\partial_xY +\bar\rho_+Y.
\end{aligned}
\end{equation}
 Furthermore, the Rankine-Hugoniot conditions
(\ref{DiffJshock}) and (\ref{shockspeedeq}) can be written as
\begin{equation}\label{1stUSRH}
Y_t=-\mA_1(Y_x, s(t)-x_0),
\end{equation}
and
\begin{equation}\label{shockspeedY}
\begin{aligned}
s'=\mA_2(Y_x,s-x_0),
\end{aligned}
\end{equation}
respectively.
Moreover,  direct computation yields
\begin{equation*}
\begin{aligned}
Y(s(t),t)&=E_+(s(t),t)-\bar E_+(s(t))\\
&=\bar E_-(s(t))-\bar E_+(s(t))\\
&=\bar E_-(x_0)+\bar E_-(s(t))-\bar E_-(x_0)-\bar E_+(x_0)-\bar
E_+(s(t))+\bar E_+(x_0)\\
&=\partial_x\bar E_-(x_0)\cdot(s(t)-x_0)-\partial_x\bar
E_+(x_0)\cdot(s(t)-x_0)+R_1.
\end{aligned}
\end{equation*}
Using the third equation in the Euler-Poisson equations (\ref{UnsteadyEP}), one has
\begin{equation}\label{shockY}
s(t)-x_0=\mA_3(Y(t, s(t))),
\end{equation}
where $\mA_3(0)=0$ and
\begin{equation*}
\frac{\partial\mA_3}{\partial Y}=\frac{1}{\bar{\rho}_-(x_0)-\brhop(x_0)}.
\end{equation*}
Combining (\ref{1stUSRH})  and (\ref{shockY}) together yields
\begin{equation}\label{BCYshock}
\begin{aligned}
\partial_tY=\mA_4(Y_x, Y),\,\,\,\,\,\,\text{at}\,\,x=s(t),
\end{aligned}
\end{equation}
where
\begin{equation}
\mA_4(0,0)=0,\,\,\,\, \frac{\partial \mA_4}{\partial
Y_x}=\frac{c^2(\bar\rho_+)(x_0)-\bar u_+^2(x_0)}{2\bar
u_+(x_0)},\,\,\,\, \frac{\partial \mA_4}{\partial Y}=-\frac{\bar
E_+(x_0)}{2\bar u_+(x_0)}.
\end{equation}
Note that on the right boundary, $x=L$, $Y$ satisfies
\begin{equation}\label{BCYR}
\partial_xY=0 ,\,\,\,\,\,\,\text{at}\,\, x=L.
\end{equation}
Our goal is to derive uniform estimates for $Y$ and $s$ which
satisfy (\ref{EQY}), (\ref{shockY}) (\ref{BCYshock}) and
(\ref{BCYR}).

Introduce the transformation
\begin{equation}
\tilde{t}=t,\,\,\,\, \tilde{x}=(L-x_0)\frac{x-s(t)}{L-s(t)}+x_0,\,\,\,\, \sigma(\tilde{t})=s(t)-x_0,
\end{equation}
to transform the problem to the fixed domain $[x_0, L]$. Let
\begin{equation}\label{q1q2}
q_1(\tilde x, \s)=\frac{L-\tilde x}{L-x_0-\sigma (\tilde t)},\,\,\,\, q_2(\s)=\frac{L-x_0}{L-x_0-\sigma (\tilde t)}.\end{equation}
It is easy to verify that
$$\frac{\partial}{\partial t}=\frac{\partial}{\partial \tilde t}-\sigma '(\tilde t)q_1\frac{\partial}{\partial{\tilde x}},\qquad  \frac{\partial}{\partial x}=q_2\frac{\partial}{\partial{\tilde x}},$$
\begin{align*}
\frac{\partial^2}{\partial t^2}&=\frac{\partial^2}{\partial \tilde t^2}+(q_1 \sigma' (\tilde t))^2\frac{\partial^2}{\partial \tilde x^2}-2q_1\sigma'(\tilde t)\frac{\partial^2}{\partial \tilde x\partial \tilde t}-q_1\left(\sigma''(\tilde t)+2\frac{( \sigma' (\tilde t))^2}{L-x_0-\sigma (\tilde t)}\right)\frac{\partial}{\partial \tilde x},\end{align*}
$$\frac{\partial^2}{\partial x\partial t}=q_2\left (\frac{\partial^2}{\partial \tilde x\partial \tilde t}+\frac{\sigma'(\tilde t)}{L-x_0-\sigma (\tilde t)}\frac{\partial}{\partial{\tilde x}}-q_1\sigma'(t)\frac{\partial^2}{\partial \tilde x^2}\right),\,\, \frac{\partial^2}{\partial \tilde x^2}=q_2^2\frac{\partial^2}{\partial \tilde x^2}.$$
So (\ref{3120}) becomes
\begin{equation}\label{eqinnew}
\begin{aligned}
&\partial_{\tilde t\tilde t}Y+q_2\partial_{\tilde x}\left(p(\bar\rho_+)+\frac{\bar J^2}{\bar\rho_+}-p(\bar\rho_++q_2 Y_{\tilde x})-\frac{(\bar J-Y_{\tilde t}-q_2 Y_{\tilde x})^2}{\bar\rho_++q_2 Y_{\tilde x}}\right)\notag\\
&-2\sigma'(t) q_1 \partial_{\tilde x\tilde t}Y+(q_1\sigma'(t))^2\partial_{\tilde x \tilde x}Y
+\bar\rho_+Y+\bar E_+q_2\partial_{\tilde x}Y+q_2 YY_{\tilde x}\notag\\
&-2\frac{(\sigma'(\tilde t))^2}{L-x_0-\sigma(\tilde t)}q_1 \partial_{\tilde x}Y\notag\\
= &\sigma''(\tilde t)q_1 \partial_{\tilde x}Y.
\end{aligned}
\end{equation}

With the help of straightforward calculations, the equation (\ref{shockY}) becomes
\begin{equation}\label{firstshockeq}
\s = \mA_3(Y(t, \tilde{x}=x_0)).
\end{equation}
The
equation for the shock front,  (\ref{shockspeedY}),  becomes
\begin{equation}
\frac{d\sigma}{d\tilde{t}}=\mA_2(q_2(\s)Y_{\tx},\ \sigma(t)).
\end{equation}
Applying (\ref{firstshockeq}) to represent the quadratic terms for
$\sigma$ in terms of $Y$, we have, at $\tilde x=x_0$,
\begin{equation}\label{shockseq}
\frac{d\sigma}{d\tilde{t}}+\frac{\bEp}{2\bar{u}_+}(x_0)\sigma=\mC_2(Y_{\tilde{x}}, Y),
\end{equation}
where $\mC_2$ satisfies
\begin{equation*}
\left| \mC_2(Y_{\tilde{x}}, Y) +
\frac{c^2(\brhop)-\bup^2}{2(\brhop-\bar{\rho}_-)\bup} (x_0)Y_{\tilde{x}}    \right| \leq C(Y_{\tilde{x}}^2+Y^2).
\end{equation*}
Clearly, it follows from (\ref{firstshockeq}) and (\ref{shockseq}) that one can represent $\s$ and $\s'$ in terms of $Y$ and its derivatives at $\tilde{x}=x_0$. Thus, after manipulating (\ref{BCYshock}) with (\ref{firstshockeq}) and (\ref{shockseq}), we have
\begin{equation}\label{defC1}
Y_{\tilde {t}}=\mC_1(Y_{\tx}, Y), {~\rm at~} \tx=x_0.
\end{equation}
Or, by the implicit function theorem again, equivalently,  we have
\begin{equation}
Y_{\tx}=\mC_3(Y_{\tilde{t}}, Y), \quad \text{at} \,\, \tx=x_0,
\end{equation}
where $\mC_3$ satisfies
\begin{equation}
\left|\mC_3(Y_{\tilde{t}}, Y)- \frac{2\bup}{c^2(\brhop)-\bup^2}(x_0)Y_{\tilde{t}}-\frac{\bEp}{c^2(\brhop)-\bup^2}(x_0)Y\right|\leq C(Y_{\tilde{t}}^2+Y^2).
\end{equation}

In the following, we drop $\tilde{  }$ in $\tilde x$ and $\tilde t$ for simplicity of notation.

The problem takes the following compact form
\begin{equation}
\left\{
\begin{aligned}
&\mL(x,{Y},\sigma){Y}
=\sigma''(t)q_1 \partial_{ x}Y, \quad (t,x)
\in [0,\infty)\times [x_0, L],\\
&\partial_{x}{Y}=d_1(Y_t, Y)Y_t+e_1(Y_t, Y) Y,\,\, \text{at} \,\, x=x_0,\\
&\partial_{x} Y=0,\,\, \text{at}\,\, x=L,\\
&\s(t)=\mA_3(Y(t, x_0)),
\end{aligned}
\right.
\end{equation}
where,  by  using $\xi_0$ and $\xi_1$ to denote $t$ and $x$, respectively,
\begin{equation}
\begin{aligned}
\mL(x,Y, \sigma)Z =&\sum_{i,j=0}^1a_{ij}(x,Y,\nabla Y,\sigma,
\sigma' )\partial_{ij}Z+\sum_{i=0}^1b_i(x,Y,\nabla Y,\sigma,
\sigma')
\partial_i Z\\
&+g(x,Y,\nabla Y, \sigma, \sigma')Z,
\end{aligned}
\end{equation}
with
\begin{equation}
\begin{aligned}
&d_1(Y_t, Y)= \int_0^1\frac{\partial \mC_3}{ \partial
Y_t} ( \theta Y_t,\theta
Y)d\theta,\,\,e_1( Y_t, Y) = \int_0^1\frac{\partial
\mC_3}{\partial Y} ( \theta
Y_{t},\theta Y)d\theta.
\end{aligned}
\end{equation}
Furthermore, one has
\begin{equation}
\begin{aligned}
&\mL(x,0,0)Z=\mL_0Z,
\end{aligned}
\end{equation}
and
\begin{equation}\label{positivity}
\begin{aligned}
&a_{00}(x,Y, \nabla Y, \sigma, \sigma')=1,\,\, a_{01}(x,0,0,0,0)=a_{10}(x,0,0,0,0)=\frac{\bar J}{\brhop}=\bar u_+,\,\,\\
 & a_{11}(x,0,0,0,0)=-\left(p'(\brhop)-\frac{\bar J^2}{\brhop^2}\right),\\
&b_0(x,0,0,0, 0)=\partial_x\left(\frac{2\bar J}{\brhop}\right), \,\, b_1(x,0,0,0, 0)=-\partial_x\left(p'(\brhop)-\frac{\bar J^2}{\brhop^2}\right)+\bEp,\\
& g(x,0,0,0, 0)=\brhop,\,\,
d_1( 0, 0)=\frac{2\bar
u_+}{c^2(\bar\rho_+)-\bar u_+^2}(x_0),\,\, e_1(0, 0)= \frac{ \bEp}{c^2(\bar\rho_+)-\bar u_+^2}(x_0).
\end{aligned}
\end{equation}

%\newpage

\subsection{Linear Estimate}\label{seclinear}

In this subsection, we  study the  linearized problem.
\begin{thm}\label{thmlinear}
Assume that $\bar{E}_+$ satisfies (\ref{stabilitycondition}).
Let $Y$ be a smooth solution of
 the linearized problem
\begin{equation}\label{lzeropb}
\left\{
\begin{aligned}
&\mL(x,0,0)Y=0,\quad  x_0<x<L, \, \, t>0,  \\
&\partial_ x Y=\frac{2\bup}{p'(\brhop)-\bup^2}(x_0)\partial_x Y+ \frac{\bEp}{p'(\brhop)-\bup^2}(x_0)Y,\,\,\,\,\,\,\text{at}\,\,x=x_0,\\
&\partial_xY=0 ,\,\,\,\,\,\,\text{at}\,\, x=L,\\
&Y(0,x)=h_1(x), \,\, Y_t(0,x)=h_2(x),\, \, x_0<x<L.
\end{aligned}
\right.
\end{equation}
Then there exist $\alpha_0\in (0,1)$  and $T>0$ such that
\begin{equation}\label{energydsp}
\hat{\varphi}_k(Y, t+T)<\alpha_0\hat{\varphi}_k(Y, t)\qquad \text{for}\,\, t\geq 0,
\end{equation}
where $\hat{\varphi}_k$ is defined as follows
\begin{equation}
\hat{\varphi}_k(Y, t)=\sum_{m=0}^{k}\varphi_m(Y, t)
\end{equation}
with
\begin{equation}
\begin{aligned}
&\varphi_m(Y, t)=\frac{\bEp}{\brhop}(x_0)(\partial_t^mY)^2(t, x_0)\\
&\quad+\int_{x_0}^L\frac{1}{\brhop}\left\{(\partial_t^{m+1}Y)^2+
\left(p'(\bar{\rho}_+)-\frac{\bar{J}_+^2}{\bar{\rho}_+^2}\right)(\partial_x\partial_t^mY)^2
+\bar{\rho}_+\partial_t^mY^2\right\}(t,x) dx.
\end{aligned}
\end{equation}
\end{thm}

\begin{rmk}
When $\bEp$ satisfies (\ref{stabilitycondition}), it follows from the Sobolev inequality that
\begin{equation}
\begin{aligned}
\varphi_m(Y, t)>C\int_{x_0}^L\left\{(\partial_t^{m+1}Y)^2+
(\partial_x\partial_t^mY)^2
+\partial_t^mY^2\right\}(t,x) dx.
\end{aligned}
\end{equation}
This is the key reason that we can handle the case that $\bEp(x_0)$ is a negative number with small magnitude.
\end{rmk}

\begin{pf} The proof of the theorem has four steps.

Step 1: Energy estimates.  Multiplying the first equation in (\ref{lzeropb}) with $\frac{1}{\brhop(x)}\partial_tY$ on both sides and integrating by parts, we have
\begin{equation}
\begin{aligned}
0=&\frac{1}{2}\partial_t\int_{x_0}^L \frac{1}{\brhop(x)}
\left\{(\partial_tY)^2+\left(p'(\bar{\rho}_+)-\frac{\bar{J}_+^2}{\bar{\rho}_+^2}\right)(\partial_xY)^2
+\bar{\rho}_+Y^2\right\}\\
&+\int_{x_0}^L\left(-\partial_x \left(\frac{1}{\brhop}\right)\frac{\bar{J}_+}{\bar{\rho}_+}
+\frac{1}{\brhop}\partial_x\left(\frac{\bar{J}_+}{\bar{\rho}_+}\right)\right) (\partial_tY)^2dx\\
&+\int_{x_0}^L\left(\partial_x \left(\frac{1}{\brhop}\right)\left(p'(\bar{\rho}_+)-\frac{\bar J^2}{\brhop^2}\right)
+\frac{1}{\brhop}\bEp\right)\partial_tY\partial_xYdx\\
&+\frac{\bar J}{\brhop^2}(\partial_tY)^2\Big|_{x_0}^L
-\frac{1}{\brhop}\left(p'(\brhop)-\frac{\bar J^2}{\brhop^2}\right)\partial_tY\partial_xY\Big |_{x_0}^L\\
=&\sum_{i=1}^4 I_i.
\end{aligned}
\end{equation}
Noting that $\bar J$ is a constant, and using the equation (\ref{SteadyEP}), we have $I_2=I_3=0$.
Apply the boundary conditions (\ref{lzeropb}) to get
\begin{equation}
\begin{aligned}
I_4=&\frac{\bar J}{\brhop^2}(\partial_tY)^2\Big|_{x=L}-
\frac{\bar J}{\brhop^2}(\partial_tY)^2\Big|_{x=x_0}\\
&+\frac{1}{\brhop}\left(p'(\brhop)-\frac{\bar J^2}{\brhop^2}\right)\partial_tY\left(\frac{2\bar{u}_+}{p'(\brhop)-J^2/\brhop^2}\partial_t
Y+\frac{\bar{E}_+}{p'(\brhop)-J^2/\brhop^2}Y\right)\Big|_{x=x_0}\\
=&\frac{\bar J}{\brhop^2}(\partial_tY)^2\Big|_{x=L}+
\frac{\bar J}{\brhop^2}(\partial_tY)^2\Big|_{x=x_0}
+\partial_t\left(\frac{\bEp(x_0)}{2\brhop(x_0)}Y^2\right)\Big|_{x=x_0}.
\end{aligned}
\end{equation}
Thus,
\begin{equation}\label{energyidentity}
\varphi_0(Y, t)+D_0(Y, t)=\varphi_0(Y, 0),
\end{equation}
where
\begin{equation}\label{defD0}
D_0(Y, t)=2\left(\int_{0}^t \frac{\bar J}{\brhop^2}(\partial_tY)^2(s, L)ds+
\int_{0}^t \frac{\bar J}{\brhop^2}(\partial_tY)^2(s, x_0)ds\right).
\end{equation}

Step 2: Rauch-Taylor type estimates. Using the boundary condition at $x=x_0$ and the fact that
$\frac{c^2(\bar\rho_+)-\bar u_+}{2u_+}(x_0)\ge C$ for some constant $C>0$, it is easy to see that
\begin{equation}\label{estimateDzero}
D_0(Y, t)\geq C_1\int_{0}^t (Y_t^2+Y_x^2)(s, x_0)ds-C_2\int_{0}^t Y^2(s, x_0)ds,
\end{equation}
for some positive constants $C_1$ and $C_2$, independent of $t$.
Therefore,
\begin{equation}
\varphi_0(Y, t)+C_1\int_{0}^t (Y_t^2+Y_x^2)(s, x_0)ds\leq \varphi_0(Y, 0)+C_2 \int_{0}^t Y^2 (s,x_0)ds
\end{equation}
Following from the argument in \cite{Rauch} and the details in Section \ref{SecRauch}, there exists a $T>0$ and $\delta\in (0, T/4)$ such that
\begin{equation}\label{estRT1}
\begin{aligned}
\int_{0}^T (Y_t^2+Y_x^2)(t, x_0)dt \geq &\int_{\frac{T}{2}-\delta}^{\frac{T}{2}+\delta} \varphi_0(Y, s)ds -C_3\int_0^T Y^2(t, x_0)dt.
\end{aligned}
\end{equation}
Since $\varphi_0(Y, t)$ is decreasing with respect to $t$, in view of (\ref{energyidentity}), we have
\begin{equation} \label{estimateRT}
\begin{aligned}
\int_{0}^T (Y_t^2+Y_x^2)(t, x_0)dt  \geq & \delta \varphi_0(Y, \frac{T}{2}+\delta)-C_3\int_0^T Y^2(t, x_0)dt\\
\geq & \delta \varphi_0(Y, T)-C_3\int_0^TY^2(t, x_0)dt.
\end{aligned}
\end{equation}
Combining (\ref{energyidentity}), (\ref{estimateDzero}) and (\ref{estimateRT}) together, we have
\begin{equation}\label{energycpt}
(1+C_4)\varphi_0(Y, T)\leq \varphi_0(Y, 0)+C_5\int_{0}^T Y^2 (t, x_0) dt,
\end{equation}
for some positive constants $C_4$ and $C_5$, independent of $t$.

Step 3: Spectrum of the solution operator. Let us define a new norm $\|\cdot\|_{\bX}$ for the function $h=(h_1,h_2) \in H^1\times L^2([x_0,L])$,
\begin{equation}
\|h\|_{\bX}^2= \frac{\bEp}{\brhop}(x_0)|h_1|^2(x_0)+\int_{x_0}^L\frac{1}{\brhop}\left\{|h_2|^2+
(p'(\bar{\rho}_+)-\frac{\bar{J}_+^2}{\bar{\rho}_+^2})|h_1'|^2
+\bar{\rho}_+|h_1|^2\right\}(x) dx.
\end{equation}
The associated complex Hilbert space will be denoted by $(\bX, \|\cdot\|_{\bX})$.
Define the solution operator $S_t: \bX\mapsto \bX$ as
\begin{equation}
S_t(h)=(Y(t,\cdot),Y_t(t,\cdot))
\end{equation}
where $Y$ is the solution of the problem (\ref{lzeropb}) with the initial data $h=(h_1,h_2)$. By  (\ref{energyidentity}), we can see that $S_t$ is bounded and satisfies
\begin{equation}
\|S_t\|\leq 1.
\end{equation}
It follows from the spectrum radius theorem \cite{Lax} that the radius of the spectrum of $S_t$ is less than or equal to 1, i.e.,  $|\mathcal{\sigma}(S_t)|\leq 1$.
Furthermore, we can define a map $K: \bX\mapsto L^2([0, T])$ as
\begin{equation}
K(h)=Y(t,x_0).
\end{equation}
Thus the estimate (\ref{energycpt}) can be written as
\begin{equation}
(1+C_6)\|S_T(h)\|_{\bX}\leq \|h\|_{\bX}+C_7\|Kh\|_{L^2([0, T])},
\end{equation} for some positive constants $C_6$ and $C_7$.
Note that for the initial data $h\in \bX$, there exists a solution $Y\in H^1([0, T]\times [x_0, L])$.
It follows from the trace theorem of $H^1$ space, \cite{Evans}, that $K$ is compact. Applying the Lemma on page 81 in \cite{Rauch2}, we have the following proposition:
\begin{prop}\label{PropWeyl}
There are only finite generalized eigenvalues for the operator $S_T$ in the annulus $\{\frac{1}{1+C_6}<|z|\leq 1\}$ on the complex plane, each of these eigenvalues has the finite multiplicity.
\end{prop}

Step 4: Refine estimate for the spectrum of the solution operator.
\begin{prop}\label{Propeigen}
There is no generalized eigenvalues of $S_T$ on the circle $|z|=1$.
\end{prop}
\noindent {\bf Proof of the proposition \ref{Propeigen}:} Suppose not, there exist $\omega\in \mathbb{R}$ and $V\in \bX$ such that
\begin{equation}
(S_T-e^{i\omega} I)V=0.
\end{equation}
Note that the identity (\ref{energyidentity}) still holds in the complex setting if we replace the square terms in $\varphi_0$ and $D_0$ by the square of modulus. Thus
\begin{equation}
\varphi_0(Y, 0)- \varphi_0(Y, nT)=n D_0(Y, T).
\end{equation}
Since both $\varphi_0(Y, nT)$ and $\varphi_0(Y, 0)$ are both positive and finite, it yields that
\begin{equation}\label{zerodsp}
D_0(Y, T)=0.
\end{equation}
Therefore
\begin{equation}
\varphi_0(Y, t)=\varphi_0(Y, 0),\,\, \text{for all}\,\,t.
\end{equation}
Let
\begin{equation}
\mBV=ker(S_T-e^{i\omega}I).
\end{equation}
Note that the coefficients in the problem (\ref{lzeropb}) do not depend on $t$, so
\[
(S_T-e^{i\omega}I)S_t=S_t(S_T-e^{i\omega}I)
\]
In particular, $S(t)\mBV\subset\mBV$ for any $t$, so $\mBV$ is invariant with respect to $S_t$.
Therefore,  $S_t|_{\mBV}$ is a semigroup on a finite dimensional subspace.
This yields that
\begin{equation}
S(t)|_{\mBV}=e^{tA}
\end{equation}
for some matrix $A$. Choose an eigenvector $w$ of $A$ such that
\begin{equation}
Aw=\lambda w.
\end{equation}
So
\begin{equation}
S_t w=e^{\lambda t}w.
\end{equation}
Note that if  $S_t$ has an eigenvalue on the unit circle, then there is a solution with $\lambda$  pure imaginary. If $S_t w=e^{i\beta t}w$, then it follows from (\ref{zerodsp}) that $w(x_0)=0$. Furthermore, by the boundary condition at $x_0$, we have $w'(x_0)=0$. Note that $w$ is an eigenfunction, so $w=0$.
This finishes the proof of the proposition.
\qed

Combining Proposition \ref{PropWeyl} and Proposition \ref{Propeigen} together, one has
\begin{equation}
\mathfrak{\sigma}(S_T)\subset\{|z|\leq \sqrt{\alpha_0}\}\,\, \text{for}\,\, 0<\alpha_0<1,
\end{equation}
where $\mathfrak{\sigma}(S_T)$ is the spectrum of $S_T$.
Since the coefficients for both the equation and the boundary conditions are independent of $t$, applying this estimate to the time derivatives of the solutions implies that
\begin{equation}
\varphi_m(Y, T)\leq \alpha_0\varphi_m(Y, 0).
\end{equation}
Summing these estimates together, we have (\ref{energydsp}).

This finishes the proof of the theorem.
\end{pf}

%\newpage

The consequence of the theorem is that the solution of the problem
(\ref{lzeropb}) decays exponentially.

\begin{cor}\label{cordecay}
There exist constants $\lambda_0>0$ and $C>0$ such that for any
solution $Y$ of the problem (\ref{lzeropb}),  one has
\begin{equation}
\hat{\varphi}_k(Y,t)\leq Ce^{-\lambda_0 t} \hat{\varphi}_k (Y, 0)
\end{equation}
and
\begin{equation}
\int_0^{\infty} e^{\frac{\lambda_0 t}{4}}\sum_{l=1}^{k+1} ( |\partial_t^l Y|^2(t, x_0)+|\partial_t^l Y|^2(t, L)) dt\leq C \hat{\varphi}_k(Y, 0).
\end{equation}
\end{cor}
\begin{pf}
Noting that $\hat{\varphi}_k(Y, t)$ is decreasing in $t$,
for any $t\in [nT, (n+1)T)$, $n\in \mbN$, one has
\begin{equation}
\begin{aligned}
\hat{\varphi}_k(Y, t)\leq & \hat{\varphi}_k (Y, nT)\leq \alpha_0^n \hat{\varphi}_k (Y, 0)\\
\leq & \alpha_0^{\frac{t}{T}-1}\hat{\varphi}_k(Y, 0) = e^{-\lambda_0 t} \alpha_0^{-1}\hat{\varphi}(Y,0),
\end{aligned}
\end{equation}
where we have chosen $\lambda_0 =-\frac{\ln \alpha_0}{T}$.

Note that
\begin{equation}
\begin{aligned}
&\int_{2^i T}^{2^{i+1}T} e^{\frac{\lambda_0 t}{4}}\sum_{l=1}^{k+1}(|\partial_t^l Y|^2 (t, x_0)  + |\partial_t^l Y|^2 (t, L)) dt\\
\leq & e^{\lambda_0 2^{i-1}T} \int_{2^i T}^{2^{i+1}T} \sum_{l=1}^{k+1}(|\partial_t^l Y|^2 (t, x_0)  + |\partial_t^l Y|^2 (t, L)) dt\\
= & e^{\lambda_0 2^{i-1}T} (\hat{\varphi}_k (Y, 2^i T) -\hat{\varphi}(Y,2^{i+1}T))\\
\leq & e^{\lambda_0 2^{i-1}T} \alpha_0^{-1} e^{-\lambda 2^i T}\hat{\varphi}_k (Y, 0)= \alpha_0^{-1} e^{-\lambda 2^{i-1} T}\hat{\varphi}_k (Y,0),
\end{aligned}
\end{equation}
where we use the energy equality (\ref{energyidentity}).
Thus
\begin{equation}
\begin{aligned}
&\int_0^{\infty}  e^{\frac{\lambda_0 t}{4}}( |\partial_t^l Y|^2(t, x_0)+|\partial_t^l Y|^2(t, L)) dt\\
= & \sum_{i=0}^{\infty}\int_{2^i T}^{2^{i+1}T} e^{\frac{\lambda_0 t}{4}}\sum_{l=1}^{k+1}(|\partial_t^l Y|^2 (t, x_0)  + |\partial_t^l Y|^2 (t, L)) dt\\
\leq & \sum_{i=0}^{\infty} \alpha_0^{-1} e^{-\lambda_0 T 2^{i-1}} \hat{\varphi}_k(Y, 0)\leq C\hat{\varphi}_k(Y, 0).
\end{aligned}
\end{equation}
This finishes the proof of the corollary.
\end{pf}

%\newpage

\subsection{Uniform A Priori Estimates}\label{secuap}

The existence of local-in-time solutions for the problem  (\ref{UnsteadyEP}), (\ref{stead}), (\ref{boundary1}), (\ref{supsub1}), (\ref{rh1}), and (\ref{vacuum}) can be achieved via the argument in \cite{LY}.
 In order to get global existence of the nonlinear  problem  (\ref{UnsteadyEP}), (\ref{stead}), (\ref{boundary1}), (\ref{supsub1}), (\ref{rh1}), and (\ref{vacuum}), we need only to derive global a priori estimate for the problem
\begin{equation}\label{nonlinear}
\left\{
\begin{aligned}
&\mL(x,Y,\sigma)Y=\sigma''(t)q_1(x, \s) \partial_x Y, \quad x_0<x<L, \,\, t>0, \\
&\partial_xY=d_1(Y_t, Y)Y_t+e_1(Y_t, Y),\,\,\,\,\,\,\text{at}\,\,x=x_0,\\
& \s(t) =\mA_3(Y(t,x_0)),\\
&\partial_xY=0 ,\,\,\,\,\,\,\text{at}\,\, x=L,\\
&Y(0,x)=h_1(x), \,\, Y_t(0,x)=h_2(x),\quad x_0<x<L.
\end{aligned}
\right.
\end{equation}

For $t>T$ in Theorem \ref{thmlinear} and an integer $k\ge 15$, we introduce the following notations
\begin{equation}
\begin{aligned}
\lbn (Y, \sigma)\rbn =&\sup_{\tau\in [0, t)}\sum_{0\leq m\leq k-6}\left(\sum_{0\leq l\leq m} e^{\frac{\lambda \tau}{16}}\|
\partial_t^l \partial_x^{m-l}Y(\tau, \cdot) \|_{L^{\infty}([x_0, L])} +e^{\frac{\lambda \tau}{16}}\left|\frac{d^m \s}{dt^m}\right|(\tau)\right)\\
&  + \sup_{0\leq \tau \leq t} \left( \sum_{\substack{ 0\leq l \leq m,\\0\leq m\leq k  } }  \|\partial_t^l \partial_x^{m-l}Y(\tau, \cdot) \|_{L^{2}([x_0, L])}   +  \sum_{ 0\leq l \leq k}   \|\partial_t^l \partial_x^{k+1-l}Y(\tau, \cdot) \|_{L^{2}([x_0, L])}  \right)\\
&  + \sup_{0\leq \tau \leq t}   \|\partial_t^{k+1}Y(\tau, \cdot)-\frac{d^{k+1}\s}{d t^{k+1}}q_1(\cdot, \s)\partial_x Y(\tau, \cdot) \|_{L^{2}([x_0, L])} \\
& +\sum_{\substack{0\leq
l\leq m,\\ 0\leq m\leq k+1}} \|\partial_t^l\partial_x^{m-l}
Y(\cdot,x_0)\|_{L^2[0, t]} +\sum_{\substack{0\leq
l\leq m,\\ 0\leq m\leq k+1}} \|\partial_t^l\partial_x^{m-l}
Y(\cdot,L)\|_{L^2[0, t]}\\
&+\sum_{0\leq m \leq k+1}\left\|\frac{d^m\sigma}{dt^m}\right\|_{L^2[0, t]},
\end{aligned}
\end{equation}
where $\lambda$ is defined in (\ref{deflam}).

In this subsection, when there are no specific indications,  we always assume $a_{ij}$,  $b_i$, and $g$ are functions of $(x, Y, \nabla Y, \sigma, \dot{\sigma})$, $d_1$ and $e_1$  are functions of $(Y, Y_t)$.
Furthermore, for any $l\in \mathbb{N}$ and given $Y$ and $\s$ such that $\lbn (Y, \s)\rbn<\infty$, we define
\begin{equation}
\hat{\Phi}_l(Z, t; Y, \sigma)=\tPhi_{l-1}(Z, t; Y, \sigma)+\Phi_0(\partial_t^l Z-q_1(x, \s) Y_x\frac{d^l \s}{dt^l}, t; Y, \s),
\end{equation}
and
\begin{equation}
\tPhi_l(Z, t; Y, \s)=\sum_{m=0}^{l}\Phi_m(Z, t; Y, \s),
\end{equation}
where
\begin{equation}
\begin{aligned}
\Phi_m(Z, t; Y, \sigma)=\frac{e_1(\partial_t^mZ)^2}{d_1\brhop}(t,x_0)+ \int_{x_0}^L \frac{(\partial_t^{m+1}Z)^2 - a_{11}(\partial_x\partial_t^{m}Z)^2 +g(\partial_t^mZ)^2}{\brhop(x)}(t,x)dx.
\end{aligned}
\end{equation}

By (\ref{positivity}), it is easy to see that
\begin{equation}
\Phi_m(Z, t; Y, \sigma)(t)\ge C \int_{x_0}^L (\partial_t^{m+1}Z)^2 + (\partial_x\partial_t^{m}Z)^2 +(\partial_t^mZ)^2)(t,x)dx,
\end{equation}
for some constant $C>0$ independent of $t$ if $\big\| \|(Y, \sigma)\|\big\|\leq \epsilon$ for some $\epsilon>0$.
Also,
\begin{equation}
\hat{\Phi}_k(Z, t; 0, 0)=\hat{\varphi}_k(Z, t).
\end{equation}

\begin{prop}\label{propapest}
Assume that $\bar{E}_+$ satisfies (\ref{stabilitycondition}). There
exists an $\epsilon_0>0$, for any $\epsilon \in (0, \epsilon_0)$,
if  $(Y, \s)$ is a smooth solution of the problem (\ref{nonlinear})
with $(h_1, h_2)$ satisfying
\begin{equation}
|\s_0|+ \|h_1\|_{H^{k+2}}+\|h_2\|_{H^{k+1}}\leq \epsilon^2\leq \epsilon_0^2
\end{equation}
and $\lbn (Y, \s)\rbn\leq \epsilon$ for $t>T$, then
\begin{equation}\label{apest}
\lbn (Y,\s)\rbn \leq \frac{\epsilon}{2}.
\end{equation}
\end{prop}
\begin{pf}
We divide the proof into four steps.

Step 1: Lower order energy estimate. Taking the $m$-th ($0\leq m \leq k-1$) order derivative for the equation (\ref{nonlinear}) with respect to
$t$, then
\begin{equation}
\begin{aligned}
{\mL}(x, Y, \sigma)\partial^m_tY=\mF_m(x,Y, \sigma)+\tilde{\mF}_m(x,Y, \sigma),
\end{aligned}
\end{equation}
where
\begin{equation}
\begin{aligned}
\mF_m(x,Y, \sigma)=
\sum_{1\leq l \leq m}C_m^l\left\{-\sum_{i,j=0}^1\partial_t^{l}a_{ij}\partial_{ij}\partial^{m-l}_tY- \sum_{i=0}^1\partial_t^{l}b_{i}\partial_{i}\partial^{m-l}_tY
- \partial_t^{l}g\partial^{m-l}_tY\right\},
\end{aligned}
\end{equation}
and
\begin{equation}
\tilde{\mF}_m(x,Y, \sigma)= \sum_{0\leq l \leq m}C_m^l\frac{d^{l+2}\sigma}{dt^{l+2}}\partial_t^{m-l}
\left(q_1(x, \s) Y_x\right),
\end{equation}
where and in the following $C_m^l$ is the binomial coefficient,
\begin{equation*}
C_m^l=\left\{
\begin{aligned}
&\frac{m!}{l!(m-l)!},\quad \text{if}\,\, m\geq l,\\
&0, \qquad \qquad\quad \text{if}\,\, m<l.
\end{aligned}
\right.
\end{equation*}
Multiplying the both sides by $\partial_t^{m+1}Y/\brhop(x)$ and
integrating on $\Omega=:[0, t]\times [x_0, L]$ for $t\in [0, T]$, noticing that $a_{00}=1$, we get
\begin{equation}
\begin{aligned}
&\iint_{\Omega} \mL(x, Z, \sigma)\partial^m_tY \frac{\partial_t^{m+1}Y}{\brhop(x)}d\tau dx\\
=&\int_{x_0}^L \frac{(\partial_t^{m+1}Y)^2 - a_{11}(\partial_x\partial_t^{m}Y)^2 +g(\partial_t^mY)^2}{(2\brhop(x))}(t,x)dx\\
&+\iint_{\Omega} \left(\frac{b_0}{\brhop(x)}-\partial_x\left(\frac{a_{01}}{\brhop(x)}\right)\right) (\partial_t^{m+1}Y)^2 d\tau dx\\
&+\iint_\Omega \left(\frac{b_1}{\brhop(x)}-\partial_x(\frac{a_{11}}{\brhop})\right)\partial_t^m\partial_xY \partial_t^{m+1}Y d\tau dx\\
&+\iint_{\Omega}\frac{\partial_ta_{11}}{\brhop}\frac{(\partial_x \partial_t^{m}Y)^2}{2}(t,x)- \frac{\partial_t g}{\brhop(x)}\frac{(\partial_t^mY)^2}{2}(\tau,x) d\tau dx\\
&+\int_0^t \left(\frac{a_{11}\partial_x\partial_t^mY\partial_t^{m+1}Y}{\brhop(x)}
+\frac{a_{01}}{\brhop(x)} (\partial_t^{m+1}Y)^2\right)(\tau, L)dt\\
&-\int_0^t\left(\frac{a_{11}\partial_x\partial_t^mY\partial_t^{m+1}Y}{\brhop(x)}
+\frac{a_{01}}{\brhop(x)} (\partial_t^{m+1}Y)^2\right)(\tau, x_0)d\tau\\
&-\int_{x_0}^L \frac{(\partial_t^{m+1}Y)^2 - a_{11}(\partial_x\partial_t^{m}Y)^2 +g(\partial_t^mY)^2}{2\brhop(x)}(0,x)dx\\
=&\sum_{i=1}^7J_i.
\end{aligned}
\end{equation}

In view of the fact
\begin{equation*}
\begin{aligned}
&\frac{b_0(x,0,0,0,0)}{\brhop(x)}-\partial_x(\frac{a_{01}(x,0,0,0,0)}{\brhop(x)})\\
=&\frac{1}{\brhop(x)}\partial_x\left(\frac{2\bar J}{\brhop}\right)-
\partial_x\left(\frac{\bar J/\brhop}{\brhop}\right)\\
=&0
\end{aligned}
\end{equation*}
and
\begin{equation*}
\begin{aligned}
&\frac{b_1(x,0,0,0,0)}{\brhop(x)}-\partial_x\left(\frac{a_{11}(x,0,0,0,0)}{\brhop(x)}\right)\\
=&\frac{1}{\brhop(x)}\left[-\partial_x\left(p'(\brhop)-\frac{\bar J^2}{\brhop^2}\right)+\bEp\right]
-\partial_x\left(\frac{-(p'(\brhop)-\bar J^2/\brhop^2)}{\brhop}\right)\\
=&0,
\end{aligned}
\end{equation*}
we conclude
\begin{equation*}
\begin{aligned}
|J_2+J_3|\leq  & C \int_0^t \left[\sum_{l=0}^2 \left( \sum_{i=0}^l \|\partial_t^i \partial_x^{l-i} Y\|_{L^{\infty}([x_0, L])}+|\partial_t^l \s|(\tau) \right)\right]  \sum_{l=1}^{m+1} \sum_{i=0}^1 \|\partial_t^{l-i} \partial_x^{i}  Y(\tau, \cdot) \|_{L^2[x_0, L]}^2 d\tau.
\end{aligned}
\end{equation*}
Note that
\begin{equation}
 \sum_{l=1}^{m+1} \sum_{i=0}^1 \|\partial_t^{l-i} \partial_x^{i}  Y(\tau, \cdot) \|_{L^2[x_0, L]}^2  \leq  C \tPhi_m(Y, \tau; Y,\s), \end{equation}
thus
\begin{equation}\label{estJ23}
\begin{aligned}
|J_2+J_3|\leq  &C  \sup_{0\leq \tau \leq t} \left| e^{\frac{\lambda \tau}{64}} \left[\sum_{l=0}^2 \left( \sum_{i=0}^l \|\partial_t^i \partial_x^{l-i} Y\|_{L^{\infty}([x_0, L])}+|\partial_t^l \s|(\tau) \right)\right] \right| \times \\
& \int_0^t   e^{-\frac{\lambda \tau}{64}}  \tPhi_m(Y, \tau; Y, \sigma) d\tau\\
\leq & C \big\| \|(Y, \s)\|\big\| \int_0^t   e^{-\frac{\lambda \tau}{64}}  \tPhi_m(Y, \tau; Y, \sigma)  d\tau
\end{aligned}
\end{equation}
Similarly,
\begin{equation}\label{estJ4}
|J_4|\leq C\big\|\|(Y, \sigma)\rbn \int_0^t e^{-\frac{\lambda \tau}{64}}  \tPhi_m(Y, \tau; Y, \sigma)  d\tau.
\end{equation}
It follows from the boundary condition at $x=L$, $Y_x(t, L)=0$, that
\begin{equation}\label{estJ5}
J_5=\int_0^t \frac{a_{01}}{\bar\rho_+}(\partial_t^{m+1}Y)^2(\tau, L)dt>0,
\end{equation}
if $\lbn (Y, \s)\rbn$ is sufficiently small.

Differentiating the boundary condition
\begin{equation}
\partial_x Y=d_1(Y_t, Y)Y_t +e_1(Y_t, Y)Y
\end{equation}
$m$ times with respect to $t$, one has
\begin{equation}\label{diff1stbc}
\partial_t^m Y_x= d_1\partial^{m+1}_tY +e_1\partial_t^m Y
+\mG_m,\,\,\text{at}\,\, x=x_0,
\end{equation}
where $\mG_m$ satisfies
\begin{equation}
|\mG_m|\leq C\left(\sum_{l=0}^{[\frac{m}{2}]+1}\sum_{i=0}^l  |\partial_t^{l} Y||\partial_t^{m+1-i}Y| + \sum_{l=0}^{[\frac{m}{2}]+1}\ |\partial_t^{l} Y|^2\right).
\end{equation}
Therefore, at $x=x_0$,
\begin{equation}
\begin{aligned}
&\frac{a_{11}\partial_x\partial_t^mY\partial_t^{m+1}Y+a_{01}
(\partial_t^{m+1}Y)^2}
{\brhop} \\
=& \frac{(a_{11}{d_1}+a_{01})(\partial_t^{m+1}Y)^2-
\partial_t(a_{11}e_1\frac{(\partial_t^m Y)^2}{2})
+\partial_t(a_{11}e_1)\frac{(\partial_t^m Y)^2}{2}}{\brhop}
+\frac{\mG_m \partial_t^{m+1}Y}{\brhop}
\end{aligned}
\end{equation}
So
\begin{equation}\label{estJ6}
\begin{aligned}
J_6=&-\int_0^t\frac{(a_{11}{d_1}+a_{01})}{\brhop}(\partial_t^{m+1}Y)^2(\tau, x_0)d\tau +\left( a_{11}e_1\frac{(\partial_t^m Y)^2}{2}\right)(t, x_0) \\
&-\int_0^t\partial_t\left(a_{11}e_1\right)\frac{(\partial_t^m Y)^2}{2\brhop}(\tau, x_0)d\tau+\int_0^t\frac{\mG_m \partial_t^{m+1}Y}{\brhop}(\tau,x_0)d\tau\\
&-\left(a_{11}e_1\frac{(\partial_t^m Y)^2}{2}\right)(t, 0)\\
\geq&\int_0^t\frac{-(a_{11}{d_1}
+a_{01})}{\brhop}(\partial_t^{m+1}Y)^2(\tau, x_0)d\tau\\
& +\left(a_{11}e_1 \frac{(\partial_t^m Y)^2}{2}\right)(t, x_0)
-\left(a_{11}e_1\frac{(\partial_t^m Y)^2}{2}\right)(0, x_0)\\
&-C\big\|\|(Y, \sigma)\|\big\| \int_0^t \sum_{l=0}^{m+1} (\partial_t^l Y)^2(\tau, x_0) d\tau.
\end{aligned}
\end{equation}

Summing up all the estimates (\ref{estJ23}), (\ref{estJ4}), (\ref{estJ5}), and (\ref{estJ6}) together yields
\begin{equation}\label{estleftL}
\begin{aligned}
&\iint_{\Omega} \hat{\mL}(x, Y, \sigma)\partial^m_tY \frac{\partial_t^{m+1}Y}{\brhop(x)}dtdx\\
\geq &\Phi_m(Y, t; Y, \sigma)-\Phi_m(Y, 0; Y, \sigma)+\mD_m(Y, t; Y, \sigma)\\
& -C\big\|\|(Y, \sigma)\|\big\|  \left[   \int_0^t e^{-\frac{\lambda \tau}{64}}
 \tPhi_m(Y, \tau; Y, \sigma) d\tau  + \tmD_m(Y, t; Y,\s)   \right],
\end{aligned}
\end{equation}
where $\tmD_m$ is defined as follows
\begin{equation}
\tmD_m=\sum_{l=0}^m \mD_l(Z, t; Y, \s)
\end{equation}
and
\begin{equation}
\begin{aligned}
\mD_l(Z, t; Y, \sigma)=\int_{0}^{t}
\frac{-(a_{11}{d_1}+a_{01})}{\brhop} (\partial_t^{l+1}Z)^2(\tau,
x_0)d\tau  + \int_{0}^{t}\frac{
a_{01}}{\brhop}(\partial_t^{l+1}Z)^2(\tau, L)d\tau.
\end{aligned}
\end{equation}
In particular, $\mD_0(Z,t; 0, 0)=D_0(Z, t)$ which is defined in (\ref{defD0}).

Note that
\begin{equation}
\begin{aligned}
& - \iint_{\Omega}\mF_m\frac{\partial_t^{m+1}Y}{\brhop}dtdx\\
=&\iint_{\Omega} \sum_{\substack{1\leq l \leq m,\\  0\le i, j\le 1}}C_m^l\frac{\partial_t^{l}a_{ij}}{\brhop(x)}
\partial_{ij}\partial^{m-l}_tY \partial_t^{m+1}Y+ \sum_{\substack{1\leq l \leq m,\\  0\le i\le 1}}C_m^l
\frac{\partial_t^{l}b_{i}}{\brhop(x)}\partial_{i}\partial^{m-l}_tY
\partial_t^{m+1}Y\\
&+\iint_{\Omega}\sum_{1\leq l\leq m}C_m^l\frac{\partial_t^{l}g}{\brhop(x)}
\partial^{m-l}_t Y \partial_t^{m+1}Y\\
=&\iint_{\Omega} \sum_{\substack{1\leq l \leq 6,\\ 0\le i, j\le 1}}C_m^l\left(\frac{\partial_t^{l}a_{ij}}{\brhop(x)}
\partial_{ij}\partial^{m-l}_tY + \frac{\partial_t^{l}b_{i}}{\brhop(x)}
\partial_{i}\partial^{m-l}_tY \right)\partial_t^{m+1}Y\\
&+\iint_{\Omega}\sum_{1\leq l\leq 6}C_m^l\frac{\partial_t^{l}g}{\brhop(x)}
\partial^{m-l}_t Y \partial_t^{m+1}Y\\
&+\iint_{\Omega} \sum_{\substack{7\leq l \leq m, \\ 0\le i, j\le 1}}C_m^l \left(\frac{\partial_t^{l}a_{ij}}{\brhop(x)}
\partial_{ij}\partial^{m-l}_tY +\frac{\partial_t^{l}b_{i}}{\brhop(x)}
\partial_{i}\partial^{m-l}_tY\right) \partial_t^{m+1}Y\\
&+\iint_{\Omega}\sum_{7\leq l\leq m}\frac{\partial_t^{l}g}{\brhop(x)}
\partial^{m-l}_t Y \partial_t^{m+1}Y\\
=&\sum_{i=1}^4 K_i.
\end{aligned}
\end{equation}
We estimate each term as follows. The first part of $K_1$ can be estimated as follows
\begin{equation}
\begin{aligned}
&\left|\iint_{\Omega} \sum_{\substack{1\leq l \leq 6,\\ 0\le i, j\le 1}}C_m^l\frac{\partial_t^{l}a_{ij}}{\brhop(x)}
\partial_{ij}\partial^{m-l}_tY \partial_t^{m+1}Y\right|\\
\leq & C \int_0^T \sum_{1\leq l \leq 6}\left\|
\frac{\partial_t^{l}a_{ij}}{\brhop(x)}\right\|_{L^{\infty}([x_0, L])}\|
\partial_{ij}\partial^{m-l}_tY \|_{L^2([x_0, L])}
\|\partial_t^{m+1}Y\|_{L^2([x_0, L])}d\tau
\end{aligned}
\end{equation}
Using the equation (\ref{nonlinear}), we have
\begin{equation}
\sum_{l=1}^{m+1}\sum_{i=0}^l \|\partial_t^i \partial_x^{l-i} Y(\tau, \cdot) \|_{L^2[x_0, L]}^2 \leq C \left(\tPhi_m(Y, \tau; Y, \s) +\lbn (Y, \s)\rbn^2  e^{-\frac{\lambda \tau}{16}} \sum_{l=1}^{m+1} \left|\frac{d^l\s}{dt^l}\right|^2  \right).
\end{equation}
Thus
\begin{equation}
\begin{aligned}
&\left|\iint_{\Omega} \sum_{\substack{1\leq l \leq 6,\\ 0\le i, j\le 1}}C_m^l\frac{\partial_t^{l}a_{ij}}{\brhop(x)}
\partial_{ij}\partial^{m-l}_tY \partial_t^{m+1}Y\right|\\
\leq& C\big\|\|(Y, \sigma)\|\big\| \int_0^t  e^{-\frac{\lambda \tau}{64}}
 \left[ \tPhi_m(Y, \tau; Y, \sigma) + \lbn (Y, \s)\rbn  \sum_{l=0}^{m+1} \left|\frac{d^l \s}{dt^l}\right|^2  \right]
d\tau.
\end{aligned}
\end{equation}
The estimates for the second part of  $K_1$ and $K_2$ are quite similar to the estimate above.
Moreover, for the first part of $K_3$, we have
\begin{equation}
\begin{aligned}
&\left|\iint_{\Omega} \sum_{\substack{7 \leq l \leq m,\\ 0\le i, j\le 1}} C_m^l\frac{\partial_t^{l}a_{ij}}{\brhop(x)}
\partial_{ij}\partial^{m-l}_tY \partial_t^{m+1}Y\right|\\
\leq & \int_0^t  \sum_{7\leq l \leq m} \left \|
\frac{\partial_t^{l}a_{ij}}{\brhop(x)}\right\|_{L^2([x_0,L])} \left\| \partial_{ij}
\partial^{m-l}_tY \right\|_{L^{\infty}([x_0,L])}
\|\partial_t^{m+1}Y\|_{L^2([x_0,L])}d\tau\\
\leq & C \big\|\|(Y, \sigma)\|\big\| \int_0^t   e^{-\frac{\lambda \tau}{64}}  \left[ \tPhi_m(Y, \tau; Y, \sigma)   + \lbn (Y, \s)\rbn  \sum_{l=0}^{m+1} \left|\frac{d^l \s}{dt^l}\right|^2  \right]
d\tau.
\end{aligned}
\end{equation}
So, we have
\begin{equation}\label{estFm}
\begin{aligned}
\iint_{\Omega}\mF_m\frac{\partial_t^{m+1}Y}{\brhop}dtdx \leq C \big\|\|(Y, \sigma)\|\big\| \int_0^t   e^{-\frac{\lambda  \tau}{64}}  \left[ \tPhi_m(Y, \tau; Y, \sigma)   + \lbn (Y, \s)\rbn  \sum_{l=0}^{m+1} \left|\frac{d^l \s}{dt^l}\right|^2  \right]
d\tau.
\end{aligned}
\end{equation}

Similarly,
\begin{equation}\label{esttFm}
\begin{aligned}
&\left|\iint_{\Omega}\tilde{\mF}_m\frac{\partial_t^{m+1}Y}{\brhop}dtdx\right|\\
\leq
& C\iint_{\Omega} \sum_{0\leq l \leq 6}  \frac{d^{l+2}\s}{dt^{l+2}}
\partial_t^{m-l}q_1(x, \s)
\frac{\partial_t^{m+1}Y}{\brhop}dtdx\\
&+C \iint_{\Omega} \sum_{7 \leq l \leq m}  \frac{d^{l+2}\s}{dt^{l+2}}
\partial_t^{m-l}q_1(x, \s)
\frac{\partial_t^{m+1}Y}{\brhop}dtdx\\
\leq & C\int_{0}^t \sum_{0\leq l \leq 6}  \left| \frac{d^{l+2}\s}{dt^{l+2}}\right|\cdot \left\|
\partial_t^{m-l}q_1(x, \s) \right\|_{L^2([x_0,L])}
\left\|\frac{\partial_t^{m+1}Y}{\brhop}\right\|_{L^2([x_0,L])}dt\\
&+ C \int_{0}^t \sum_{7 \leq l \leq m} \left| \frac{d^{l+2}\s}{dt^{l+2}}\right|\cdot
\left\|\partial_t^{m-l}q_1(x, \s)\right\|_{L^{\infty}([x_0,L])}
\left\|\frac{\partial_t^{m+1}Y}{\brhop}\right\|_{L^2([x_0,L])}dt\\
\leq &C \lbn (Y, \s)\rbn \int_0^t  e^{-\frac{\lambda \tau}{64}}
 \left[ \tPhi_m(Y, \tau; Y, \sigma) + \lbn (Y, \s)\rbn  \sum_{l=0}^{m+2} \left|\frac{d^l \s}{dt^l}\right|^2  \right]
 d\tau\\
 & +  C \lbn (Y, \s)\rbn \int_0^t  e^{-\frac{\lambda \tau}{64}}
 \left[ \tPhi_m(Y, \tau; Y, \sigma) +  \sum_{l=0}^{m+2} \left|\frac{d^l \s}{dt^l}\right|^2  \right]
 d\tau\\
 \leq  &  C \lbn (Y, \s)\rbn \int_0^t  e^{-\frac{\lambda \tau}{64}}
 \left[ \tPhi_m(Y, \tau; Y, \sigma) +  \sum_{l=0}^{m+2} \left|\frac{d^l \s}{dt^l}\right|^2  \right]
 d\tau,
\end{aligned}
\end{equation}
where we have used the estimate
\begin{equation*}
\left\|\frac{d^{l+2}\s}{dt^{l+2}}\right\|_{L^2([x_0, L])}\leq C \left|\frac{d^{l+2}\s}{dt^{l+2}}\right|.
\end{equation*}

Combining the estimates (\ref{estleftL}),  (\ref{estFm}) and (\ref{esttFm}) together, one has
\begin{equation}\label{nllenergy}
\begin{aligned}
&\tPhi_m(Y, t; Y, \sigma)+ \tmD_m(Y, t; Y, \sigma)\leq \tPhi_m(Y, 0; Y, \sigma) \\
&\quad  +  C \big\|\|(Y, \sigma)\|\big\| \left[ \int_0^t
 e^{-\frac{\lambda \tau}{64}}
 \left[ \tPhi_m(Y, \tau; Y, \sigma) +   \sum_{l=0}^{m+2} \left|\frac{d^l \s}{dt^l}\right|^2  \right]  d\tau +\tmD_m(Y, t; Y, \s)\right],
\end{aligned}
\end{equation}
for $m=0, 1, \cdots ,k-1$.

Step 2: The highest order energy estimates. Take the $k$-th order derivative for the equation (\ref{nonlinear}) with respect to
$t$, then
\begin{equation}
\begin{aligned}
&\mL(x, Y, \sigma)\partial^k_tY =\mF_k(x,Y, \sigma)\\
&+ \frac{d^{k+2}\s}{dt^{k+2}}q_1(x, \s)Y_x +
\sum_{0\leq l \leq k-1}C_k^l\frac{d^{l+2}\s}{dt^{l+2}}\partial_t^{k-l}
\left(q_1(x, \s)Y_x\right).
\end{aligned}
\end{equation}
In order to handle the term $\frac{d^{k+2}\sigma}{dt^{k+2}}$, we rewrite the above equation as
\begin{equation}\label{cptnlheq}
{\mL}(x, Y,\nabla Y, \sigma,\dot{\sigma})\check{Y}=\mF_k(x,Y, \sigma)+\check{\mF}(x,Y, \sigma),
\end{equation}
where
\begin{equation}
\check{Y}=\partial_t^k Y-q_1(x, \s) Y_x \frac{d^k \s}{dt^k}
\end{equation}
and
\begin{equation}
\begin{aligned}
&\check{\mF}(x,Y, \sigma)
=\sum_{0\leq l \leq k-1}C_k^l\frac{d^{l+2}\s}{dt^{l+2}}\partial_t^{k-l}
(q_1(x, \s)Y_x)\\
&-\left(\frac{d^{k+1}\s}{dt^{k+1}}
\partial_t\left(q_1(x, \s)Y_x\right) +\frac{d^{k}\s}{dt^{k}}
\partial_t^2\left(q_1(x, \s)Y_x\right)\right)\\
&-(2a_{01}\partial_t\partial_x+a_{11}\partial_x^2
+\sum_{i=0}^1b_i\partial_i +g) \left( \frac{d^{k}\s}{dt^{k}}
q_1(x, \s)Y_x\right).
\end{aligned}
\end{equation}

Note that the boundary condition for $\check{Y}$  at $x=x_0$ is
\begin{equation}
\begin{aligned}
\partial_x \check{Y} = &(d_1+\frac{\partial d_1}{\partial Y_t}Y_t +\frac{\partial e_1}{\partial Y_t} Y)\partial_t \check{Y} +e_1\check{Y}\\
& + (d_1+\frac{\partial d_1}{\partial Y_t}Y_t +\frac{\partial e_1}{\partial Y_t} Y)q_1 Y_x \frac{d^{k+1}\s}{dt^{k+1}}\\
& + (d_1+\frac{\partial d_1}{\partial Y_t}Y_t +\frac{\partial e_1}{\partial Y_t} Y)\partial_t(q_1 Y_x) \frac{d^{k}\s}{dt^{k}} +e_1 q_z Y_x \frac{d^k \s}{dt^k} -\partial_x (q_1 Y_x) \frac{d^k \s}{dt^k}\\
& +\sum_{l=1}^{k-1} C_{k-1}^l ( \partial_t^l(\frac{\partial d_1}{\partial Y_t}) Y_t + \partial_t^l(\frac{\partial e_1}{\partial Y_t}) Y) \partial_t^{k+1-l}Y\\
&+\partial_t^{k-1}(\frac{\partial d_1}{\partial Y} \partial_tY)\partial_tY  +\partial_t^{k-1}(\frac{\partial e_1}{\partial Y} \partial_tY)Y\\
& +\sum_{l=1}^{k-1}C_k^l \partial_t^ld_1(Y_t, Y) \partial_t^{k-l+1}Y +\sum_{l=0}^{k-1}C_k^l \partial_t^l e_1(Y_t, Y)\partial_t^{k-l} Y.
\end{aligned}
\end{equation}
The associated boundary condition for $\check{Y}$ at  $x=L$ is
\begin{equation}
\partial_x \check{Y} =-\frac{\partial }{\partial x}\left(\frac{L-x}{L-x_0-\s(t)}Y_x \frac{d^k \s}{dt^k}\right)=0.
\end{equation}

Multiplying the both sides of (\ref{cptnlheq}) by $\frac{\partial_t\check{Y}}{\brhop}$ and
integrating on $\Omega$, and noting that
\begin{equation}
\sum_{l=0}^{k+1} \|\partial_t^l\partial_x^{k-l}Y(\tau, \cdot)\|^2_{L^2[x_0, L]}\leq C\left(\hPhi_k (Y, \tau; Y, \s) +\lbn (Y, \s)\rbn^2 \sum_{l=0}^{k+1}\left|\frac{d^l \s}{dt^l}\right|^2 \right),
\end{equation}
one has
\begin{equation}\label{nlhenergy}
\begin{aligned}
&\Phi_0(\check{Y}, t; Y, \sigma)+ \mD_0(\check{Y}, t; Y, \sigma)\leq \Phi_0(\check{Y}, 0; Y, \sigma)\\
&+C\big\| \|(Y, \sigma)\|\big\| \left[
\int_0^t e^{-\frac{\lambda \tau}{64}} \left[\hat\Phi_k(Y, \tau; Z, \sigma) +\sum_{l=0}^{k+1} \left|\frac{d^{l}\s}{dt^{l}}\right|^2 \right]d\tau +\hmD_k(Y, t; Y, \s)\right],
\end{aligned}
\end{equation}
where $\hmD_k$ is defined as follows
\begin{equation}
\hat{\mD}_k(Z, t; Y, \s)=\tmD_{k-1}(Z, t; Y, \s) +\mD_0(\partial_t^k Z-q_1(x, \s)Y_x\frac{d^k \s}{dt^k}, t; Y, \s).
\end{equation}
Adding  the estimates (\ref{nllenergy}) and (\ref{nlhenergy}) yields
\begin{equation}\label{nlfullenergy}
\begin{aligned}
&\hat{\Phi}_k(Y, t; Y, \sigma)+ \hat{\mD}_k(Y, t; Y, \sigma)\leq \hat{\Phi}_k(Y, 0; Y, \sigma)\\
&+C\big\| \|(Y, \sigma)\|\big\| \left\{
\int_0^t  e^{-\frac{\lambda \tau}{64}} \left[\hat\Phi_k(Y, \tau; Y, \sigma) +\sum_{l=0}^{k+1}\left|\frac{d^{l}\s}{dt^{l}}\right|^2 \right] d\tau +\hmD_k (Y, t; Y, \s)\right\}.
\end{aligned}
\end{equation}

%\newpage

Step 3: Boundedness  of the energy. Differentiating  the equation for the shock front
\begin{equation}\label{eqshockfront}
\s(t) =\mA_3(Y(t, x_0))
\end{equation}
with respect to $t$,  we have
\begin{equation}\label{estsspeed}
\sum_{l=0}^{k+1} \left|\frac{d^l\s}{dt^l} \right|^2 (\tau) \leq C\left( |\s(\tau)|^2 +\sum_{l=0}^{k+1}
|\partial_t^l Y(\tau, x_0)|^2   \right).
\end{equation}
Therefore, we have
\begin{equation}\label{estsig}
\sum_{l=0}^{k+1} \left|\frac{d^l\s}{dt^l} \right|^2 \leq C
\sum_{l=1}^{k} |\partial_t^l Y(\tau, x_0)|^2 +|\partial_t \check{Y}|^2.
\end{equation}
Combining (\ref{estsspeed}) and (\ref{estsig}) together, we get
\begin{equation}
\int_0^t e^{-\frac{\lambda \tau}{64}} \sum_{l=0}^{k+1}
\left|\frac{d^{l}\s}{dt^{l}}\right|^2  d\tau \leq C \left(\int_0^t
e^{-\frac{\lambda \tau}{64}}\Phi_0(Y, \tau; Y, \s) d\tau +\hat{\mD}_k(Y,t;
Y, \s)\right)
\end{equation}
Thus the energy estimate (\ref{nlfullenergy})  is equivalent to
\begin{equation}
\begin{aligned}
&\hat{\Phi}_k(Y, t; Y, \sigma)+ \hat{\mD}_k(Y, t; Y, \sigma)\leq  \hat{\Phi}_k(Y, 0; Y, \sigma)\\
&+C\big\| \|(Y, \sigma)\|\big\|\left ( \hat{\mD}_k(Y, \tau; Y,
\sigma ) + \int_0^t  e^{-\frac{\lambda \tau}{64}}
\hat\Phi_k(Y, \tau; Y, \sigma) d\tau   \right).
\end{aligned}
\end{equation}
If $\big\|\|(Y, \s)\|\big\|\leq  \epsilon$, then
\begin{equation}
\hat{\Phi}_k(Y, t; Z, \sigma)+  \hat{\mD}_k(Y, t; Z, \sigma)\leq 2 \hat{\Phi}_k(Y, 0; Z, \sigma)\leq C\epsilon^4.
\end{equation}
This yields
\begin{equation}\label{energybound}
\begin{aligned}
&\sup_{0\leq \tau \leq t} \left( \sum_{\substack{ 0\leq l \leq m,\\0\leq m\leq k  } }  \|\partial_t^l \partial_x^{k-l}Y(\tau, \cdot) \|_{L^{2}([x_0, L])}   +  \sum_{ 0\leq l \leq k}   \|\partial_t^l \partial_x^{k+1-l}Y(\tau, \cdot) \|_{L^{2}([x_0, L])}  \right)\\
&  + \sup_{0\leq \tau \leq t}   \|\partial_t^{k+1}Y(\tau, \cdot)-\frac{d^{k+1}\s}{d t^{k+1}}q_1(\s)\partial_x Y(\tau, \cdot) \|_{L^{2}([x_0, L])} \\
& +\sum_{ 1\leq l\leq k+1} \|\partial_t^l
Y(\cdot,x_0)\|_{L^2[0, t]} +\sum_{1\leq l\leq k+1} \|\partial_t^l\partial_x^{m-l}
Y(\cdot,L)\|_{L^2[0, t]}\\
&+\sum_{1 \leq m \leq k+1}\left\|\frac{d^m\sigma}{dt^m}\right\|_{L^2[0, t]}\\
\leq & C\left(\sup_{0\leq \tau\leq t}\hPhi^{1/2}(Y, \tau; Y, \s) +\hmD_k^{1/2}(Y, t; Y, \s) \right)\\
\leq & C\epsilon^{2}\leq \frac{\epsilon}{4}\,\,.
\end{aligned}
\end{equation}

%\newpage

Step  4: Decay of the  lower energy and  the  shock position.  The basic idea to get the decay is to control the deviation of the solution $Y$ of the problem (\ref{nonlinear}) with the solution $\bar{Y}$ of the problem (\ref{lzeropb}) via $\bY$. The contraction of the energy for $\bY$ will also yields the contraction of the energy for $Y$. This gives the decay of the lower energy of $Y$. The decay of the shock position is a consequence of the governing equation for the shock front and the decay of the lower energy.

At time $\tau=t_0$, we can choose $\bar{h}_1\in H^{k}$ and $\bar{h}_2\in H^{k-1}$ such that there exists a solution $\bY\in C^{k-1-i}([t_0,\infty); H^{i}([x_0, L]))$ of the linear problem (\ref{lzeropb}) satisfying $\bY(t_0, \cdot) =\bar{h}_1$ and $\bY_t(t_0, \cdot) =\bar{h}_2$, and $\bY$ satisfies
\begin{equation}\label{devIC1}
\sum_{l=0}^{k-1}\sum_{i=0}^l \|\partial_t^i \partial_x^{l-i}\bY(t_0, \cdot)\|_{L^2[x_0, L]} \leq C\lbn (Y, \s)\rbn
\end{equation}
for some uniform constant $C$, and
\begin{equation}\label{devIC2}
\hPhi_{k-4}(Y-\bY, t_0;Y, \s) \leq C\lbn (Y, \s)\rbn \hPhi_{k-4}(Y, t_0; Y, \s).
\end{equation}

Note that $Y-\bY$ satisfies the equation
\begin{equation}\label{eqdeviation}
\begin{aligned}
& \sum_{i,j=0}^1 a_{ij}(x, Y,\s) \partial_{ij}(Y-\bY) +\sum_{j=0}^1 b_j(x, Y,\s)\partial_j(Y-\bY)
+g(x, Y,\s) (Y-\bY) \\
= & \sum_{i,j=0}^1 (a_{ij}(x, Y,\s)-a_{ij}(x, 0, 0)) \partial_{ij}\bY+\sum_{j=0}^1 (b_j(x, Y,\s)-b_j(x, 0,0))
\partial_j \bY\\
& +( g(x, Y,\s)-g(x, 0, 0))\bY +\s''(t) q_1(x,\sigma) \partial_x Y,
\end{aligned}
\end{equation}
and the boundary conditions
\begin{equation}
\begin{aligned}
&\partial_x(Y-\bY)=d_1(Y_t, Y)\partial_t(Y-\bY) +e_1(Y_t, Y)(Y-\bY) \\
& +(d_1(Y_t, Y)-d_1(0, 0))\partial_t\bY +(e_1(Y_t, Y)-e_1(0, 0))\bY,\,\, \text{at} \,\, x=x_0,
\end{aligned}
\end{equation}
and
\begin{equation}
\partial_x(Y-\bY)=0,\,\,\,\, \text{at}\,\, x=L.
\end{equation}

Taking the $m$-th ($m\leq k-4$) order derivatives with respect to $t$ for (\ref{eqdeviation}) and multiplying the both sides of the resulting equation with $\frac{\partial_t^{m+1}(Y-\bY)}{\brhop}$ (if $m<k-4$) or $\frac{\partial_t Z}{\brhop}$ ($Z=\partial_t^{k-4}(Y-\bY)-\frac{d^{k-4}\s}{dt^{k-4}}q_1(x, \s) Y_x$  if $m=k-4$), and integrating it on $[t_0, t_0+T]\times [x_0, L]$ where $T>0$ is the one appearing in Theorem \ref{thmlinear},  the new terms other than those appearing in the energy estimates in Steps 1 and 2 should be induced by the first three terms on the right hand side in (\ref{eqdeviation}). Here we only give the highest order energy estimates.  Let us estimate the first one as follows
\begin{equation}
\begin{aligned}
&\left|\int_{t_0}^{t_0+T} \int_{x_0}^L \sum_{l=0}^{k-4} C_{k-4}^l \partial_t^l (a_{ij}(x, Y, \s)-a_{ij}(x, 0, 0)) \partial_t^{k-4-l} \partial_{ij} \bY \frac{\partial_t Z}{\brhop} dx d\tau \right|\\
\leq & C  \int_{t_0}^{t_0+T} \sum_{l=0}^{k-4} \|\partial_t^l (a_{ij}(\cdot, Y, \s)-a_{ij}(\cdot, 0, 0))\|_{L^2[x_0, L]}\times \\
&\quad  \|\partial_t^{k-4-l} \partial_{ij} \bY \|_{L^\infty[x_0, L]}
 \left\| \frac{\partial_t Z}{\brhop}\right\|_{L^2[x_0,L]}  d\tau
 \end{aligned}
 \end{equation}
Using Cauchy-Schwarz inequality, we have
\begin{equation}
\begin{aligned}
&\left|\int_{t_0}^{t_0+T} \int_{x_0}^L \sum_{l=0}^{k-4} C_{k-4}^l \partial_t^l (a_{ij}(x, Y, \s)-a_{ij}(x, 0, 0)) \partial_t^{k-4-l} \partial_{ij} \bY \frac{\partial_t Z}{\brhop} dx d\tau \right|\\
\leq & C\int_{t_0}^{t_0+T} \left[\hPhi_{k-4}(Y, \tau; Y, \s) +\sum_{l=0}^{k-3} \left|\frac{d^l \s}{dt^l}\right|^2 \right]^{1/2} \times\\
& \quad \tPhi_{k-2}^{1/2}(\bY, \tau; Y, \s) \hPhi_{k-4}^{1/2}(Y-\bY, \tau; Y, \s) d\tau\\
\leq & C  \int_{t_0}^{t_0+T} \hPhi^{1/2}_{k-4}(Y, \tau; Y, \s) \tPhi^{1/2}_{k-2}(\bY, \tau; Y, \s) \hPhi_{k-4}^{1/2}(Y-\bY, \tau; Y, \s) d\tau \\
&+C\int_{t_0}^{t_0+T} \left[\sum_{l=0}^{k-3} \left|\frac{d^l \s}{dt^l}\right|^2\right]^{1/2} \tPhi_{k-2}^{1/2}(\bY, \tau; Y, \s) \hPhi_{k-4}^{1/2}(Y-\bY, \tau; Y, \s) d\tau\\
\leq & C  \int_{t_0}^{t_0+T}  \tPhi_{k-2}^{1/2}(\bY, \tau; Y, \s)\hPhi_{k-4}^{1/2}(Y, \tau; Y, \s) \hPhi_{k-4}^{1/2}(Y-\bY, \tau; Y, \s) d\tau \\
& +C  \int_{t_0}^{t_0+T} \tPhi_{k-2}^{1/2}(\bY, \tau; Y, \s) \hPhi_{k-4}(Y-\bY, \tau; Y, \s) d\tau\\
 &+C\int_{t_0}^{t_0+T} \sum_{l=0}^{k-3} \left|\frac{d^l \s}{dt^l}\right|^2 \tPhi_{k-2}^{1/2}(\bY, \tau; Y, \s) d\tau.
\end{aligned}
\end{equation}
Similarly, we can estimate all other terms.
In summary, we have the estimate as follows:
\begin{equation}
\begin{aligned}
&\hPhi_{k-4}(Y-\bY, t_0+T; Y,\s) +\hmD_{k-4}(Y-\bY, t_0+T; Y, \s) - \hmD_{k-4}(Y-\bY, t_0; Y, \s) \\
\leq & \hPhi_{k-4}(Y-\bY, t_0; Y, \s) +C \lbn (Y, \s)\rbn   \int_{t_0}^{t_0+T}  \sum _{l=0}^{k-3} \left|\frac{d^{l}\s}{dt^l }\right|^2 d\tau \\
& +C \lbn (Y, \s)\rbn  \int_{t_0}^{t_0+T}   \sum_{i=0}^{k-3} \|\partial_t^i \partial_x^{k-3-i} (Y-\bY)\|_{L^2}^2 d\tau\\
& + C\lbn (Y, \s) \rbn \int_{t_0}^{t_0+T} \hPhi_{k-4}^{1/2}(Y, \tau; Y, \s) \hPhi_{k-4}^{1/2}(Y-\bY, \tau; Y, \s)d\tau\\
& +C \lbn (Y, \s)\rbn   \int_{t_0}^{t_0+T} \hPhi_{k-4}(Y-\bY, \tau; Y, \s)   d\tau\\
& + C\lbn (Y, \s)\rbn ( \hmD_{k-4}(Y-\bY, t_0+T; Y, \s) - \hmD_{k-4}(Y-\bY, t_0; Y, \s))\\
& + C\lbn (Y, \s)\rbn ( \hmD_{k-4}(Y, t_0+T; Y, \s) - \hmD_{k-4}(Y, t_0; Y, \s)).
\end{aligned}
\end{equation}
Using the equation (\ref{eqdeviation}), we have
\begin{equation}
\begin{aligned}
\sum_{l=1}^{k-3} \sum_{i=0}^l \|\partial_t^i \partial_x^{l-i} (Y-\bY)\|_{L^2}^2  \leq &  C  \hPhi_{k-4}(Y-\bY, \tau; Y,\s) (1+ \sum_{l=0}^{k-3}\sum_{i=0}^l \|\partial_t^i \partial_x^{l-i}\bY\|_{L^2}^2) \\
&+ C\lbn (Y, \s)\rbn \sum_{l=0}^{k-3}\left|\frac{d^l\s}{dt^l}\right|^2.
\end{aligned}
\end{equation}

By the definitions of $\Phi_0$ and $\mD_0$, there exists a $C(\alpha_0)>0$ such that
\begin{equation}
\begin{aligned}
&\frac{33+31\alpha_0}{64}\hPhi_{k-4}(Y, t_0+T; Y,\s)\\
\leq  &\tPhi_{k-4}(\bY, t_0+T; Y,\s)   +C(\alpha_0) \hPhi_{k-4}(Y-\bY, t_0+T; Y, \s)
\end{aligned}
\end{equation}
and
\begin{equation}
\begin{aligned}
& \frac{1}{128}(\hmD_{k-4}(Y, t_0+T; Y,\s)-   \hmD_{k-4}(Y, t_0; Y,\s) )\\
\leq   &  \frac{1}{64}(\tmD_{k-4}(\bY, t_0+T; Y,\s)-   \tmD_{k-4}(\bY, t_0; Y,\s) )\\
& +  C(\hmD_{k-4}(Y-\bY, t_0+T; Y,\s)-   \hmD_{k-4}(Y-\bY, t_0; Y,\s) ),
\end{aligned}
\end{equation}
where $\alpha_0\in (0, 1)$ is the constant in (\ref{energydsp}). Therefore,
\begin{equation}\label{estYk41}
\begin{aligned}
&\frac{33+31\alpha_0}{64} \hPhi_{k-4}(Y, t_0+T; Y, \s) +\frac{1}{128}(\hmD_{k-4}(Y, t_0+T; Y, \s) -\hmD_{k-4}(Y, t_0; Y,\s) ) \\
\leq & \tPhi_{k-4}(\bY, t_0+T; Y, \s) +C(\alpha_0) \hPhi_{k-4}(Y-\bY, t_0+T; Y, \s)\\
& +\frac{1}{64}( \tmD_{k-4}(\bY, t_0+T; Y, \s) - \tmD_{k-4}(\bY, t_0; Y, \s))\\
& +C( \hmD_{k-4}(Y-\bY, t_0+T; Y, \s) - \hmD_{k-4}(Y-\bY, t_0; Y, \s))
\end{aligned}
\end{equation}
Note that the Theorem \ref{thmlinear} yields
\begin{equation}\label{tPhik4}
\begin{aligned}
&\tPhi_{k-4}(\bY, t_0+T; Y, \s) +\frac{1}{64}( \tmD_{k-4}(\bY, t_0+T; Y, \s) - \tmD_{k-4}(\bY, t_0; Y, \s))\\
\leq & \frac{1+31\alpha_0}{32} \tPhi_{k-4}(\bY, t_0; Y, \s).
\end{aligned}
\end{equation}
Combining (\ref{estYk41}) and (\ref{tPhik4}) together yields
\begin{equation}
\begin{aligned}
&\frac{33+31\alpha_0}{64} \hPhi_{k-4}(Y, t_0+T; Y, \s) +\frac{1}{128}(\hmD_{k-4}(Y, t_0+T; Y, \s) -\hmD_{k-4}(Y, t_0; Y,\s) ) \\
\leq & \frac{1+31\alpha_0}{32}  \tPhi_{k-4}(\bY, t_0; Y, \s) + C(\alpha_0)\hPhi_{k-4}(Y-\bY, t_0+T; Y, \s)  \\
& +C( \hmD_{k-4}(Y-\bY, t_0+T; Y, \s) - \hmD_{k-4}(Y-\bY, t_0; Y, \s))
\end{aligned}
\end{equation}
Using (\ref{estsig}), (\ref{devIC1}), and (\ref{devIC2}), we have
\begin{equation}\label{estYk4}
\begin{aligned}
&\frac{33+31\alpha_0}{64} \hPhi_{k-4}(Y, t_0+T; Y, \s) +\frac{1}{128}(\hmD_{k-4}(Y, t_0+T; Y, \s) -\hmD_{k-4}(Y, t_0; Y,\s) ) \\
\leq &  \frac{1+31\alpha_0}{32}  \tPhi_{k-4}(\bY, t_0; Y, \s) + C  \hPhi_{k-4}(Y-\bY, t_0; Y, \s)  \\
& +C \lbn (Y, \s)\rbn   \int_{t_0}^{t_0+T}  \sum _{l=0}^{k-3} \left|\frac{d^{l}\s}{dt^l }\right|^2 d\tau\\
&+ C\lbn (Y, \s) \rbn \int_{t_0}^{t_0+T}\hPhi_{k-4}^{1/2}(Y, \tau; Y, \s) \hPhi_{k-4}^{1/2}(Y-\bY, \tau; Y, \s)d\tau\\
& +C \lbn (Y, \s)\rbn   \int_{t_0}^{t_0+T} \hPhi_{k-4}(Y-\bY, \tau; Y, \s)   d\tau\\
& + C\lbn (Y, \s)\rbn ( \hmD_{k-4}(Y, t_0+T; Y, \s) - \hmD_{k-4}(Y, t_0; Y, \s))\\
& +C\lbn (Y, \s)\rbn ( \hmD_{k-4}(Y-\bY, t_0+T; Y, \s) - \hmD_{k-4}(Y-\bY, t_0; Y, \s))\\
\leq &  \frac{1+31\alpha_0}{32}  \tPhi_{k-4}(\bY, t_0; Y, \s) + C\epsilon  \hPhi_{k-4}(Y, t_0; Y, \s)  \\
&+ C\epsilon \tPhi_{k-4}(\bY, t_0; Y, \s) \\
& +C\epsilon  \int_{t_0}^{t_0+T} \hPhi_{k-4}(Y, \tau; Y, \s)   d\tau\\
& + C \epsilon ( \hmD_{k-4}(Y, t_0+T; Y, \s) - \hmD_{k-4}(Y, t_0; Y, \s)).
\end{aligned}
\end{equation}

In the same way as (\ref{nlfullenergy}), we have the following estimate
\begin{equation}\label{fullk4}
\begin{aligned}
&\hPhi_{k-4}(Y, t; Y,\s) +(1-C\epsilon)(\hmD_{k-4}(Y, t; Y, \s) -\hmD_{k-4}(Y, t_0; Y, \s)) \\
\leq & 2 \hPhi_{k-4}(Y, t_0; Y, \s).
\end{aligned}
\end{equation}

Substituting (\ref{fullk4}) into (\ref{estYk4}) and using (\ref{devIC1}), and (\ref{devIC2}), one has
\begin{equation}\label{declowenergy}
\begin{aligned}
\frac{34+30\alpha_0}{64} \hPhi_{k-4}(Y, t_0+T; Y, \s)
\leq  \frac{2+30\alpha_0}{32}  \hPhi_{k-4}(Y, t_0; Y, \s),
\end{aligned}
\end{equation}
if $\epsilon$ is sufficiently small.
As same as the proof for Corollary  \ref{cordecay},  it follows from (\ref{declowenergy}) that one has
\begin{equation}
\hPhi_{k-4}(Y_t, t; Y, s)+\s^2(t) \leq C (\tPhi_{k-4}(Y, 0; Y ,\s)+\s^2(0)) e^{-2\lambda t},
\end{equation}
where
\begin{equation}\label{deflam}
\lambda =-\frac{\ln \frac{1+\alpha}{2}}{2T}, \,\, \, \, \text{and}\,\,\,\, \alpha= \frac{17+15\alpha_0}{2+30\alpha_0}.
\end{equation}
 Thus
\begin{equation}\label{maxnormdecay}
\sum_{l=0}^{k-6}\|Y(t,\cdot)\|_{L^{\infty}[x_0, L]}\leq C\epsilon^2 e^{-\lambda t},
\end{equation}
if we have $\hPhi_{k-4}(Y, 0; Y,\s) \leq \epsilon^4$.
This yields that
\begin{equation}\label{shockdecay}
\sum_{l=0}^{k-6}\left|\frac{d^l \s}{dt^l}\right|\leq  C\epsilon^2 e^{-\lambda t}.
\end{equation}
Combining (\ref{energybound}), (\ref{maxnormdecay}), and (\ref{shockdecay}), one has (\ref{apest}).
This finishes the proof of the Proposition \ref{propapest}.
\end{pf}

Once one has the Proposition \ref{propapest}, Theorem \ref{Thmdystability} follows from the standard continuation argument and local existence result.

%\newpage

\section{Linear Instability}\label{secinsta}
We consider the linear instability for the steady transonic shock
solutions when $\bar E(x_0)<-C$ for some positive constant $C$, where $x_0\in [0, L)$ is the shock location for the steady transonic shock solution.
We consider the linearized problem (\ref{lzeropb}), which can be written in the following form:
\begin{equation}\label{lzeropb1}
\left\{\begin{aligned}
&\partial_{tt}Y-\partial_x\left(\left(c^2(\bar\rho_+)-\bup^2\right)\partial_xY\right)+\partial_x\left(2\bup\partial_tY\right)+\bar\rho_+Y+\bar E_+\partial_xY=0, \\
&\qquad \text{for}\,\,  t>0,\,\,  0<x<L, \\
&\partial_t Y=\frac{c^2(\bar \rho_+)-\bar u_+^2}{2\bar u_+}\partial_x Y-\frac{\bar E_+}{2\bar u_+}
Y,\,\, \text{at}\,\, x=x_0,\\
&\partial_x Y=0,\,\, {\rm at}\,\, x=L,
\end{aligned}
\right.
\end{equation}
where $c^2(\bar\rho_+)=p'(\bar\rho_+)$, $\bar u_+=\frac{\bar J}{\bar\rho_+}$.
Suppose that
\begin{equation}\label{instability}
\bar E_+(x_0)<-C
\end{equation}
for some positive constant $C$.
We look for solutions for the problem (\ref{lzeropb1}) of the form $Y=e^{\lambda t}Z$. Then
\begin{equation}\label{180}
\left\{
\begin{aligned}
&(p'(\bar\rho_+)-\bar u_+^2)\partial_x^2 Z+(\partial_x(c^2(\bar\rho_+)-\bar u_+^2)-2\bar u_+\lambda
-\bar E_+)Z_x \\
& \qquad -(\lambda^2+2\lambda \partial_x\bar u_++\rho) Z=0, \quad \text{for}\,\,    x_0<x<L,\\
&\partial_x Z=\frac{2\bar u_+}{c^2(\bar\rho_+)-\bar u_+^2}(x_0)(\frac{\bar E_+(x_0)}{2\bar u_+}+\lambda)Z, \text{at}\,\, x=x_0,\\
&\partial_x Z=0 ,\,\, \text{at}\,\, x=L.
\end{aligned}
\right.
\end{equation}

We use the shooting method. Given a parameter $Z(x_0)=\alpha>0$, we consider the following initial value problem:
\begin{equation}\label{181}
\left\{
\begin{aligned}
&(p'(\bar\rho_+)-\bar u_+^2)\partial_x^2 Z+(\partial_x(p'(\bar\rho_+)-\bar u_+^2)-2\bar u_+\lambda
-\bar E_+)Z_x \\
&\qquad-(\lambda^2+2\lambda \partial_x\bar u_++\rho) Z=0, \quad \text{for}\,\, x>x_0,\\
&Z(x_0)=\alpha>0, \partial_x Z=\frac{2\bar u_+}{p'(\rho_+)-\bar u_+^2}(x_0)\left(\frac{\bar E_+(x_0)}{2\bar u_+}+\lambda\right)Z, \quad \text{at}\,\, x=x_0.\\
\end{aligned}
\right.
\end{equation}

Since $\bup>0$ and $c^2(\brhop)>\bar u_+$. Therefore, in view of (\ref{instability}),  if $\lambda=0$, then $Z_x(x_0)<0$. Therefore, there exists $L_1>x_0$ such that
$Z_x(x)<0$ for $x_0\leq x\leq L_1$.

If $\lambda=-\frac{\bEp}{\bar u_+}(x_0)$, then $Z_x(x_0)>0$. then there exists $L_2>x_0$ such that
$Z_x(x)>0$ for $x_0\leq x \leq L_2$.   Let
$$L=\min\{L_1, L_2\}. $$

By the continuous dependence of ODE with respect to the initial data
and the parameters, there exists a $0<\lambda<-2\frac{E}{2\bar u_+}(x_0)$ such that the problem (\ref{181}) admits a solution $Z$ satisfying $Z_x(L)=0$ which is a solution of (\ref{180}) on $[x_0, L]$.  This implies that
the linearized problem can have exponentially growing solutions.

This finishes the proof of Theorem \ref{Thminstability}.

%\newpage

\section{Appendix: Rauch-Taylor argument}\label{SecRauch}
In this appendix, we give a detailed calculation for the estimate (\ref{estRT1}), which is motivated by the estimate in \cite{Rauch}.
 For simplicity, we introduce the following characteristic variables for the equation in  (\ref{lzeropb}):
 \begin{equation}
 \begin{aligned}
& \theta=t-\frac{1}{2}\int_{x_0}^x\frac{1}{\bup+c(\brhop)}dx-\frac{1}{2} \int_{x_0}^x\frac{1}{\bup-c(\brhop)}dx,\\
 &\zeta= \frac{1}{2}\int_{x_0}^x\frac{1}{\bup+c(\brhop)}dx-\frac{1}{2} \int_{x_0}^x\frac{1}{\bup-c(\brhop)}dx.
\end{aligned}
\end{equation}
Obviously, $\zeta=0$ if $x=x_0$. We also denote $\zeta=\zeta_L$ if $x=L$.
Then the equation in (\ref{lzeropb}) will be
\begin{equation}\label{chareq}
\partial_{\theta\theta}Y-\partial_{\zeta\zeta}Y+M\partial_{\zeta}Y+NY=0,
\end{equation}
where
\begin{equation}
M=\frac{(p'(\brhop)-\bup^2)(2p'(\brhop)-p''(\brhop)\brhop)}{2c^3(\brhop)\brhop}\frac{d\brhop}{dx},\,\, N=\frac{p'(\brhop)-\bup^2}{p'(\brhop)}\brhop.
\end{equation}
Let $Y=e^{k\zeta}Z$, then one has
\begin{equation}\label{newchareq}
\partial_{\zeta\zeta}Z-\partial_{\theta\theta}Z +(2k-M)\partial_{\zeta}Z+(k^2-kM-N)Z=0.
\end{equation}
If we Choose $k$ sufficiently large, then
\begin{equation}
2k-M>0\,\,\text{and}\,\,  k^2-kM-N>0    \,\, \text{for}\,\, \zeta\in [0,\zeta_L].
\end{equation}
Multiplying the equation (\ref{newchareq}) with $\partial_{\zeta}Z$ and integrating on the domain $\Omega=\{(\zeta,\theta)| \theta_0+\zeta\leq\theta \leq \theta_1-\zeta,  0\leq \zeta\leq \zeta_0\}$ with $\zeta_0\in [0,\zeta_L]$, then we have
\begin{equation}
\begin{aligned}
&\int_{\partial\Omega}\frac{ (\partial_{\zeta}Z)^2 n_{\zeta}-2\partial_{\zeta}Z\partial_{\theta}Z n_{\theta}+(\partial_{\theta}Z)^2 n_{\zeta}}{2} ds  +\int_{\partial\Omega}\frac{(k^2-kM-N)Z^2}{2}n_{\zeta} ds\\
&+\iint_{\Omega}(2k-M)(\partial_{\zeta}Z)^2+\frac{\partial_{\zeta}(kM+N)}{2}Z^2 d\theta d\zeta=0,
\end{aligned}
\end{equation}
where $(n_{\theta},n_{\zeta})$ is the unit outer normal on the boundary $\partial\Omega$.
Therefore, we have
\begin{equation}
\begin{aligned}
&\int_{\theta_0+\zeta_0}^{\theta_1-\zeta_0}\frac{ (\partial_{\zeta}Z)^2 +(\partial_{\theta}Z)^2 }{2} (\theta, \zeta_0)  +\frac{(k^2-kM-N)Z^2}{2}(\theta,\zeta_0) d\theta\\
&+\int_{0}^{\zeta_0}\int_{\theta_0+\zeta}^{\theta_1-\zeta}(2k-M)(\partial_{\zeta}Z)^2(\theta,\zeta)+\frac{\partial_{\zeta}(kM+N)}{2}Z^2 (\theta, \zeta)d\theta d\zeta\\
\leq & \int_{\theta_0}^{\theta_1}\frac{ (\partial_{\zeta}Z)^2 +(\partial_{\theta}Z)^2 }{2} (\theta, 0)  +\frac{(k^2-kM-N)Z^2}{2}(\theta,0) d\theta,
\end{aligned}
\end{equation}
here we require $\theta_1-\theta_0\geq 2\zeta_L$.
Integrating this inequality with respect to $\zeta_0$ on $[0, \zeta_L]$, then we have
\begin{equation}\label{RTenergy}
\begin{aligned}
&\int_0^{\zeta_L}d\zeta \int_{\theta_0+\zeta_0}^{\theta_1-\zeta_0}\frac{ (\partial_{\zeta}Z)^2 +(\partial_{\theta}Z)^2 }{2} (\theta, \zeta_0)  +\frac{(k^2-kM-N)Z^2}{2}(\theta,\zeta_0) d\theta\\
&+\int_0^{\zeta_L}d\zeta_0 \int_{0}^{\zeta_0}\int_{\theta_0+\zeta}^{\theta_1-\zeta}(2k-M)(\partial_{\zeta}Z)^2(\theta,\zeta)+\frac{\partial_{\zeta}(kM+N)}{2}Z^2 (\theta, \zeta)d\theta d\zeta\\
\leq & \zeta_L\int_{\theta_0}^{\theta_1}\frac{ (\partial_{\zeta}Z)^2 +(\partial_{\theta}Z)^2 }{2} (\theta, 0)  +\frac{(k^2-kM-N)Z^2}{2}(\theta,0) d\theta.
\end{aligned}
\end{equation}
If $k$ is sufficiently large, then
\begin{equation*}
\left|\int_0^{\zeta_L}d\zeta_0 \int_{0}^{\zeta_0}\int_{\theta_0+\zeta}^{\theta_1-\zeta}\frac{\partial_{\zeta}(kM+N)}{2}Z^2 (\theta, \zeta)d\theta d\zeta\right|  \leq \int_0^{\zeta_L}d\zeta \int_{\theta_0+\zeta_0}^{\theta_1-\zeta_0} \frac{(k^2-kM-N)Z^2}{4}(\theta,\zeta) d\theta.
\end{equation*}
Thus the estimate (\ref{RTenergy}) yields
\begin{equation*}
\begin{aligned}
&\int_0^{\zeta_L}d\zeta \int_{\theta_0+\zeta_0}^{\theta_1-\zeta_0}\frac{ (\partial_{\zeta}Z)^2 +(\partial_{\theta}Z)^2 }{2} (\theta, \zeta_0)  +\frac{(k^2-kM-N)Z^2}{4}(\theta,\zeta_0) d\theta\\
\leq & \zeta_L\int_{\theta_0}^{\theta_1}\frac{ (\partial_{\zeta}Z)^2 +(\partial_{\theta}Z)^2 }{2} (\theta, 0)  +\frac{(k^2-kM-N)Z^2}{2}(\theta,0) d\theta.
\end{aligned}
\end{equation*}
Transforming back to $x$ and $t$ coordinates, one has
\begin{equation}
\int_{x_0}^L dx\int_{t_0}^{t_1}[(\partial_t Y)^2 +(\partial_x Y)^2 +Y^2](t, x) dt\leq C\int_0^T [ (\partial_t Y)^2 +(\partial_x Y)^2 +Y^2] (t, x_0)dt
\end{equation}
for some $0<t_0<t_1<T$. Note that
\begin{equation}
\varphi_0(Y, t)\leq C\int_{x_0}^L [(\partial_t Y)^2 +(\partial_x Y)^2 +Y^2](t, x) dx,
\end{equation}
one has
\begin{equation}
\int_0^T [ (\partial_t Y)^2 +(\partial_x Y)^2 ] (t, x_0)dt\geq C\int_{\frac{T}{2}-\delta}^{\frac{T}{2}+\delta}\varphi_0(Y, t)dt -C \int_0^TY^2(t, x_0)dt
\end{equation}
for $T>0$ suitably large and  $\delta \in (0, T/4)$. This is exactly the estimate (\ref{estRT1}).

 {\bf Acknowledgments:} This work was initiated when Luo visited
The Institute of Mathematical Sciences, The Chinese University of
Hong Kong and Xie was a postdoctoral fellow there. Both of them
thank the institute's support and hospitality. Luo's research is partially supported by
an NSF grant DMS-0742834 (continuing as  DMS-0839864), Rauch's research is partially supported by
an NSF grant DMS-0807600, and Xin's research is partially supported by Hong Kong RGC Earmarked
Research Grants CUHK 4040/06P, CUHK 4042/08P, and the RGC Central Allocation Grant
CA05/06.SC01.


\bigskip

\begin{thebibliography}{99}


\bibitem{MarkPhase}
Uri~M. Ascher, Peter A. Markowich, Paola Pietra, and Christian
Schmeiser, {\it A phase plane analysis of transonic solutions for
the hydrodynamic semiconductor model}, Math. Models Methods Appl.
Sci. {\bf 1} (1991), no. 3, 347--376.

\bibitem{Shu}
D. P. Chen, R. S. Eisenberg, J. W. Jerome, C. W. Shu, {\it A
hydrodynamic model of temperature change in open ionic channels},
Biophys J. {\bf 69} (1995), 2304-2322.


\bibitem{CF1}
Gui-Qiang Chen and Mikhail Feldman,  {\it Multidimensional transonic
shocks and free boundary problems for nonlinear equations of mixed
type}, J. Amer. Math. Soc. {\bf 16} (2003), no. 3, 461--494.



\bibitem{ChenWang}
G.-Q. Chen and  D. Wang, {\it Convergence of shock capturing schemes
for the compressible {E}uler-{P}oisson equations}, Comm. Math. Phys.
{\bf 179} (1996),  no. 2, 333--364.

\bibitem{DeMark1d}
P. Degond and  P. A. Markowich, {\it On a one-dimensional
steady-state hydrodynamic model for semiconductors}, Appl. Math.
Lett. {\bf 3} (1990),  no. 3, 25--29.


\bibitem{DeMark3d}
P. Degond and  P. A. Markowich, {\it A steady state potential flow
model for semiconductors}, Ann. Mat. Pura Appl. (4)  {\bf 165}
(1993), 87--98.


\bibitem{EmbidGM}
Pedro Embid, Jonathan Goodman,  and Andrew Majda,{\it Multiple
steady states for {$1$}-{D} transonic flow}, SIAM J. Sci. Statist.
Comput.  {\bf 5} (1984), no. 1, 21--41.


\bibitem{Evans}
Lawrence C. Evans,  {\it Partial differential equations},  Graduate Studies in Mathematics, 19. American Mathematical Society, Providence, RI, 1998.


\bibitem{Gamba1d}
Irene Mart{\'{\i}}nez Gamba, {\it Stationary transonic solutions of
a one-dimensional hydrodynamic model for semiconductors}, Comm.
Partial Differential Equations {\bf 17} (1992), no. 3-4, 553--577.


\bibitem{GambaMorawetz}
Irene M. Gamba and Cathleen S. Morawetz, {\it A viscous
approximation for a {$2$}-{D} steady semiconductor or transonic gas
dynamic flow: existence theorem for potential flow}, Comm. Pure
Appl. Math. {\bf 49} (1996),  no. 10, 999--1049.

\bibitem{Ha}
Seung-Yeal Ha,  {\it $L\sp 1$ stability for systems of conservation laws with a nonresonant moving source},  SIAM J. Math. Anal. {\bf 33} (2001), no. 2, 411--439.

\bibitem{HaYang}
Seung-Yeal Ha and Tong Yang,  {\it $L\sp 1$ stability for systems of hyperbolic conservation laws with a resonant moving source}, SIAM J. Math. Anal.  {\bf 34} (2003), no. 5, 1226--1251

\bibitem{Huang}
Feimin Huang, Ronghua Pan and Huimin Yu,  {\it Large time behavior of Euler-Poisson system for semiconductor},
Sci. China Ser. A {\bf 51} (2008), no. 5, 965--972.

\bibitem{Kato}
Tosio Kato,  {\it Perturbation theory for linear operators}, Reprint of the 1980 edition, Classics in Mathematics, Springer-Verlag, Berlin, 1995.

\bibitem{Lax}
Peter D. Lax,  {\it Functional analysis}, Pure and Applied Mathematics (New York), Wiley-Interscience [John Wiley \& Sons], New York, 2002.


\bibitem{LiMark}
Hailiang Li,   Peter Markowich and Ming Mei,  {\it Asymptotic behavior of subsonic entropy solutions of the isentropic Euler-Poisson equations},  Quart. Appl. Math.  {\bf 60} (2002), no. 4, 773--796.

\bibitem{LXY1}
Jun Li, Zhouping Xin and Huicheng Yin, {\it On transonic shocks in a
nozzle with variable end pressures},  Comm. Math. Phys. {\bf 291} (2009), no. 1, 111--150.


\bibitem{LXY2}
Jun Li, Zhouping Xin and Huicheng Yin, {\it A free boundary value
problem for the Euler system and 2-D transonic shock in a large
variable nozzle},   Math. Res. Lett.  {\bf 16} (2009), no. 5, 777--796.



\bibitem{LY}
Ta Tsien Li and Wen Ci Yu, {\it Boundary value problems for quasilinear hyperbolic systems}, Duke University Mathematics Series, V. Duke University, Mathematics Department, Durham, NC, 1985.

\bibitem{Lien}
Wen-Ching Lien, {\it Hyperbolic conservation laws with a moving source}, Comm. Pure Appl. Math.  {\bf 52} (1999), no. 9, 1075--1098.

\bibitem{LiuARMA}
Tai Ping Liu, {\it Transonic gas flow in a duct of varying area},  Arch. Rational Mech. Anal. {\bf 80} (1982), no. 1, 1--18.

\bibitem{LiuTP}
Tai Ping Liu, {\it Nonlinear stability and instability of transonic
flows through a nozzle}, Comm. Math. Phys. {\bf 83} (1982), no. 2,
243--260.

\bibitem{LiuJMP}
Tai Ping Liu,  {\it Nonlinear resonance for quasilinear hyperbolic equation}, J. Math. Phys.  {\bf 28} (1987), no. 11, 2593--2602.


\bibitem{LuoNX}
Tao Luo, Roberto Natalini and Zhouping Xin, {\it Large time behavior of the solutions to a hydrodynamic model for semiconductors},  SIAM J. Appl. Math. {\bf 59} (1999), no. 3, 810--830.

\bibitem{LuoXin}
Tao Luo and Zhouping Xin, {\it Transonic shock solutions for a
system of Euler-Poisson equations}, preprint, 2009.

\bibitem{Majda}
Andrew Majda,  {\it The existence of multidimensional shock fronts},  Mem. Amer. Math. Soc. {\bf 43} (1983), no. 281.


\bibitem{MarkZAMP}
Peter A. Markowich, {\it On steady state {E}uler-{P}oisson models
for semiconductors}, Z. Angew. Math. Phys. {\bf 42} (1991), no. 3,
389--407.


\bibitem{MarkRSbook}
P. A. Markowich, C. A. Ringhofer, and C. Schmeiser, {\it
Semiconductor equations}, Springer-Verlag, Vienna, 1990.

\bibitem{Metivier}
Guy M\'{e}tivier, {\it Stability of multidimensional shocks}, Advances in the theory of shock waves, 25--103, Progr. Nonlinear Differential Equations Appl., 47, BirkhŠuser Boston, Boston, MA, 2001.

\bibitem{Pao}
C. V. Pao,  {\it Nonlinear parabolic and elliptic equations},  Plenum Press, New York, 1992.

\bibitem{Peng}
Yue-Jun Peng and Ingrid Violet,  {\it Example of supersonic
solutions to a steady state {E}uler-{P}oisson system}, Appl. Math.
Lett. {\bf 19} (2006),  no. 12, 1335--1340.

\bibitem{Rauch2}
Jeffrey Rauch,  {\it Qualitative behavior of dissipative wave equations on bounded domains},  Arch. Rational Mech. Anal.  {\bf 62} (1976),  no. 1,  77--85.


\bibitem{RauchM}
Jeffrey Rauch and Frank Massey, {\it Differentiability of solutions to hyperbolic initial-boundary value problems}, Trans. Amer. Math. Soc. {\bf 189} (1974), 303--318.


\bibitem{Rauch}
Jeffrey Rauch and Michael Taylor,  {\it Exponential decay of solutions to hyperbolic equations in bounded domains}, Indiana Univ. Math. J. {\bf 24} (1974), 79--86.


\bibitem{RXX}
Jeffrey Rauch, Chunjing Xie and Zhouping Xin, {\it Global stability of transonic shock solutions in quasi-one-dimensional nozzles}, preprint, 2010.

\bibitem{RosiniStability}
Massimiliano D. Rosini, {\it Stability of transonic strong shock
waves for the one-dimensional hydrodynamic model for
semiconductors}, J. Differential Equations {\bf 199} (2004),  no. 2,
326--351.


\bibitem{RosiniPhase}
Massimiliano D. Rosini, {\it A phase analysis of transonic solutions
for the hydrodynamic semiconductor model}, Quart. Appl. Math. {\bf
63} (2005),  no. 2, 251--268.

\bibitem{WangChen}
Dehua Wang and Gui-Qiang Chen, {\it Formation of singularities in
compressible {E}uler-{P}oisson fluids with heat diffusion and
damping relaxation}, J. Differential Equations  {\bf 144} (1998), no.
1, 44--65.

\bibitem{XinYinCPAM}
Zhouping Xin and Huicheng Yin, {\it Transonic shock in a nozzle. I.
Two-dimensional case}, Comm. Pure Appl. Math.  {\bf 58} (2005), no.
8, 999--1050.

\bibitem{XinYinJDE}
 Zhouping Xin and Huicheng Yin, {\it The transonic shock in a nozzle, 2-D and 3-D complete Euler systems},
  J. Differential Equations {\bf 245} (2008), no. 4, 1014--1085.

\bibitem{ZhangBo}
Bo Zhang,  {\it Convergence of the {G}odunov scheme for a simplified
one-dimensional hydrodynamic model for semiconductor devices}, Comm.
Math. Phys. {\bf 157} (1993), no. 1, 1--22.

\end{thebibliography}
\end{document}